\documentclass[11pt]{article}
\usepackage[T1]{fontenc}
\usepackage[latin9]{inputenc}
\usepackage[top=3cm, bottom=3cm, left=1.5cm, right=1.5cm]{geometry}
\usepackage{amsmath}
\usepackage{amsthm}
\usepackage{amssymb}
\usepackage[unicode=true,pdfusetitle,
bookmarks=true,bookmarksnumbered=true,bookmarksopen=true,bookmarksopenlevel=3,
breaklinks=false,pdfborder={0 0 1},backref=false,colorlinks=false]
{hyperref}

\usepackage[english]{babel}
\usepackage{verbatim,url}            
\usepackage{booktabs}       
\usepackage{amsfonts}       
\usepackage{nicefrac}       
\usepackage[mathscr]{eucal}
\usepackage[pdftex]{graphicx}
\usepackage[colorinlistoftodos,bordercolor=orange,backgroundcolor=orange!20,linecolor=orange,textsize=scriptsize]{todonotes}
\usepackage{subfig}
\usepackage{color, makecell, afterpage}%

\makeatletter
\theoremstyle{plain}
\newtheorem{thm}{\protect\theoremname}
\theoremstyle{plain}
\newtheorem{lem}[thm]{\protect\lemmaname}
\theoremstyle{plain}
\newtheorem{prop}[thm]{\protect\propositionname}
\newtheorem{Corollary}{Corollary}[thm]
\theoremstyle{definition}
\newtheorem{definition}{Definition}
\newtheorem*{remark}{Remark}
\usepackage{multirow}

\usepackage{arydshln}

\usepackage{lmodern}
\newcommand{\1}{\mbox{1\hspace{-1mm}I}}
\numberwithin{equation}{section}

\usepackage{calc}
\makeatother

\usepackage[round]{natbib}
\usepackage{cleveref}
\usepackage[short]{optidef}

\newcommand{\bE}{\mathbb{E}}
\newcommand{\bH}{\mathbb{H}}

\newcommand{\NN}{\mathbb{N}}

\newcommand{\RR}{\mathbb{R}}

\newcommand{\tb}{\textbf} 
\newcommand{\R}{\RR} 
\newcommand{\N}{\NN} 
\renewcommand{\b}{\mathbf}
\DeclareMathOperator*{\Id}{Id}
\newcommand{\cH}{{\mathcal H}}
\newcommand{\cL}{{\mathcal L}}
\newcommand{\cR}{{\mathcal R}}
\newcommand{\cC}{{\mathcal C}}
\newcommand{\cS}{{\mathcal S}}
\newcommand{\cN}{{\mathcal N}}

\newcommand{\bla}{\bar{\lambda}}
\newcommand{\Hk}{\mathcal{H}_K}

\newcommand{\T}{\mathscr{T}}
\newcommand{\scK}{\mathscr{K}}

\newcommand{\Sx}{\mathcal{S}^x_{[t_0,T]}}
\newcommand{\Sla}{\mathcal{S}^{\lambda}_{[t_0,T]}}

\providecommand{\norh}[1]{\lVert{#1}\rVert}
\providecommand{\noridx}[2]{\lVert{#1}\rVert_{#2}}

\providecommand{\scalh}[2]{\left\langle{#1},{#2}\right\rangle}

\providecommand{\scalidx}[3]{\left\langle{#1},{#2}\right\rangle_{#3}}

\DeclareMathOperator*{\argmin}{argmin}

\DeclareMathOperator*{\Img}{Im}
\DeclareMathOperator*{\Ker}{Ker}

\DeclareMathOperator*{\proj}{proj}

\addto\captionsamerican{\renewcommand{\lemmaname}{Lemma}}
\addto\captionsamerican{\renewcommand{\propositionname}{Proposition}}
\addto\captionsamerican{\renewcommand{\theoremname}{Theorem}}
\addto\captionsenglish{\renewcommand{\lemmaname}{Lemma}}
\addto\captionsenglish{\renewcommand{\propositionname}{Proposition}}
\addto\captionsenglish{\renewcommand{\theoremname}{Theorem}}
\providecommand{\lemmaname}{Lemma}
\providecommand{\propositionname}{Proposition}
\providecommand{\theoremname}{Theorem}

\begin{document}
	\global\long\def\1{\mbox{1\hspace{-1mm}I}}%
	
	\title{The reproducing kernel Hilbert spaces underlying linear SDE Estimation, Kalman filtering and their relation to optimal control\\
	}
	\author{Pierre-Cyril Aubin-Frankowski\footnote{pierre-cyril.aubin@inria.fr}\\
		INRIA-Département d'Informatique de l'École Normale Supérieure,\\
		PSL, Research University, Paris, France \\
		Alain Bensoussan\footnote{axb046100@utdallas.edu}\thanks{also with the School of Data Science, City University Hong Kong. Research
			supported by the National Science Foundation under grants NSF-DMS-1905449 and 2204795,
			and grant from the SAR Hong Kong RGC GRF 14301321.}\\
		International Center for Decision and Risk Analysis\\
		Jindal School of Management, University of Texas at Dallas\\
	}
	
	\maketitle
	
	\emph{In honor of Prof.\ Wendell Fleming,}
	\begin{abstract}
		It is often said that control and estimation problems are in duality. Recently, in \citet{aubin2020hard_control}, we found new reproducing kernels in Linear-Quadratic optimal control by focusing on the Hilbert space of controlled trajectories, allowing for a convenient handling of state constraints and meeting points. We now extend this viewpoint to estimation problems where it is known that kernels are the covariances of stochastic processes. Here, the Markovian Gaussian processes stem from the linear stochastic differential equations describing the continuous-time dynamics and observations. Taking extensive care to require minimal invertibility requirements on the operators, we give novel explicit formulas for these covariances. We also determine their reproducing kernel Hilbert spaces, stressing the symmetries between a space of forward-time trajectories and a space of backward-time information vectors. The two spaces play an analogue role for filtering to Sobolev spaces in variational analysis, and allow to recover the Kalman estimate through a direct variational argument. For comparison, we then recover the Kalman filter and smoother formulas through more classical arguments based on the innovation process. Extension to discrete-time observations or infinite-dimensional state, tough technical, would be straightforward.
	\end{abstract}
	
	
	
	\section{INTRODUCTION}
	\tb{Context.} In preceding papers \citep{aubin2020hard_control,aubin2020Riccati,aubin2022operator}, connections have been made between the theory of reproducing kernel Hilbert spaces and control theory. Framing control problems as optimizing over Banach vector spaces of functions with a dynamic constraint was already discussed in \citet[p255]{luenberger1968optimization}. However, in the linear-quadratic case, more specific Hilbertian structures actually emerge, the quadratic cost defining a norm over the vector space of linearly controlled trajectories. The latter has then an explicit reproducing kernel which is completely determined by the cost and dynamics. Conversely, this kernel summarizes all the information about the system, and was shown in some cases to coincide with the inverse of the solution of the backward-time differential Riccati equation \citep{aubin2020Riccati} or with the controllability Gramian \citep{aubin2020hard_control}. Note that this approach does not rely on the value function, which does not live in a vector space but rather in a max-plus one, and has thus ties with tropical kernels \citep{aubin2022tropical}.
	
	Another very related field is the problem of estimation. Its ``duality'' with control theory has been studied extensively since the seminal paper by \citet{Kalman1961NewRI}, but the word dual is unfortunately often abusively used, and does not always correspond to rigorous dual spaces in a mathematical sense or to dual min-max/max-min problems. We refer to \citet[Chapter 15]{Kailath2000} for a formalization of this duality in discrete-time and to \citet[Chapter 3]{Kim_2022phd} for valuable insights in continuous-time. Unlike in control theory, kernels have been known in estimation problems since the origin of the field. Indeed reproducing kernels are the covariances of stochastic processes in this context \citep{Parzen_1961,berlinet04reproducing}, and, in the Gaussian case, summarize all the information. In this Bayesian field of Gaussian process regression, kernels are ubiquitous, forming the deterministic counterpart of the stochastic processes \citep{Rasmussen2005-pa,kanagawa2018}. This regression framework was in particular applied to stochastic differential equations \citep[see][for an introduction and recent review]{Srkk2019} in close relation with Kalman filtering. Nevertheless, to the best of our knowledge, the reproducing kernel Hilbert spaces (RKHSs) associated with the continuous-time Gaussian processes stemming from linear stochastic differential equations (SDEs) have not been written explicitly.\\ 
	
	\tb{Main results.} To achieve the identification of the RKHSs, we extensively use the more recent theory of operator-valued kernels \citep[see][and references therein]{carmeli10vector}. We prove that the linear SDEs define two vector-valued RKHSs over the time interval: one, forward, for the reconstructible trajectories, and one, backward, for the information vectors. The values of the trajectory kernel coincide with the covariance of the minimal estimation error, and thus with the Kalman filter and Rauch-Tung-Striebel Smoother. For the information kernel, the connection is with the Gramian of observability. We provide new closed-form expressions for the covariance kernels of these Markovian Gaussian processes. By drawing upon the connection between covariance and RKHSs, we also generalize the kernel formulas obtained previously for optimal control and define two dual deterministic optimization problems associated with the smoothing task. Since we work in an estimation context, extensive care was taken to provide minimal invertibility requirements on the matrices involved. To simplify a little, we limit ourselves to finite-dimensional systems, but generalization to infinite-dimensional systems (i.e.\ nonstationnary spatiotemporal Gaussian processes, \citet{Sarkka_Solin_Hartikainen_2013,Lindgren2022}) can be done in the spirit of \citet{aubin2022operator}.\\
	
	\tb{Related work.} For discrete-time differential equations, the benefits and limitations of kernel regression, defined over the time axis, were discussed in \citet{Steinke2008}. This contrasts with the wider use of kernels, defined over observation and state spaces, as off-the-shelf tools for estimation as in \citet{Kanagawa2016-sq}. Our setting is closer to the former, with all our kernels defined over time. Characterization of Markovian Gaussian processes by the form of their kernel was recalled in \citet[Chapter 3.2]{Neveu1968ProcessusAG} and \citet[Example 2, p58]{berlinet04reproducing}, and further studied in \citet{Eisenbaum_Kaspi_2006}. The formulas for the covariance kernel of a linear SDE over the state variable but without an observation process were given in \citet[Section 6.4, p.89]{Srkk2019}. Up to our knowledge, while every linear SDE leads to a Markovian Gaussian process, the converse is not known in the general case. Concerning duality, \citet{Kim_2022phd} summarizes the traditional input-output viewpoint, with the duality appearing through the adjoint of the map sending the initial condition to the output. \citet{Kim_2022phd} then delves into a promising duality framework for the nonlinear setting. We adopt a different viewpoint here, kernel-based, focusing on covariances and Hilbert spaces of trajectories.\\
	
	The paper is structured as follows. The filtering and smoothing settings and their linear estimators are presented in Section~\ref{sec:background} and related to a general perturbed two-point boundary value problem. Section~\ref{sec:riccati-kernels} gives the formulas of the Green kernels for this problem, based on impulse-response. In Section~\ref{sec:rkhs}, these kernels are shown to be the reproducing kernels of two functional spaces, of controlled trajectories and of information vectors respectively, and two dual optimization problems defined. In Section~\ref{sec:filter}, for completeness, we rederive various known formulas for the Kalman filter and RTS smoothing. The results are finally summarized in Section~\ref{sec:summary} and compared to those previously obtained by the authors in optimal control. We provide in the Appendix implementable formulas to compute the kernels based on the Hamiltonian matrix.
	
	\section{KALMAN FILTERING AND RTS SMOOTHING IN CONTINUOUS TIME}\label{sec:background}
	\subsection{STATEMENT OF THE PROBLEM}
	We follow the presentation of \citet[Chapter 7]{Bensoussan2018}. We consider on a probability space $(\Omega,\mathcal{A},P)$
	a filtration $\mathcal{F}^{t}$ and two independent standard Wiener spaces $w(t)$,
	$b(t)$ with values in $\R^{p}$ and $\R^{m}$, and respective covariances $Q$ and $R$ in the sense that the following correlation formulas hold
	\begin{equation}
		\bE [w(s)w(t)^{*}]=\int_{t_0}^{\min(s,t)}Q(\tau)d\tau,\: \bE[b(s)b(t)^{*}]=\int_{t_0}^{\min(s,t)}R(\tau)d\tau\label{eq:2-100}
	\end{equation}
	where $w(t)^*\in \R^{p,*}$ denotes the transpose and in which $Q(\cdot)\in L^1([t_0,T],\mathcal{L}(\R^{d,*};\R^{d}))$, $R(\cdot)\in L^2([t_0,T],\mathcal{L}(\R^{m,*};\R^{m}))$ and, for some $r>0$, $R(\cdot)\succcurlyeq r\Id$ a.e. The processes are assumed to be adapted to the filtration $\mathcal{F}^{t}$. We observe a dynamic system with partial information. The dynamic system is characterized by its state $x(t)\in \R^{n}$, which evolves according to the model 
	\begin{equation}
		dx(t)=(F(t)x(t)+f(t))dt+G(t)dw(t), \quad x(t_0)=x_{0}+\xi, \quad \xi \sim \mathcal{N}(0,\Pi_{0})\label{eq:2-1}
	\end{equation}
	in which $F(\cdot)\in L^1([t_0,T],\mathcal{L}(\R^{n};\R^{n}))$, $G(\cdot)\in L^2([t_0,T],\mathcal{L}(\R^{d};\R^{n}))$ and $f(\cdot)\in L^1([t_0,T],\R^{n})$ are fixed deterministic functions, which are known, and $x_{0}\in \R^{n}$, is a
	deterministic known vector. The random variable $\xi$ has values in $\R^{n}$,
	and is assumed Gaussian with zero mean and covariance matrix $\Pi_{0}$ such that $\xi$, $w(\cdot)$, $b(\cdot)$ are mutually independent. Note that we require $R$ to be invertible but do not require it for $Q$ or $\Pi_{0}$. The state $x(t)$ is not observed, we have
	instead a continuous-time observation process $y(t)$ with values in $\R^{m}$, related to $x(t)$ by the following relation\footnote{A discrete-time observation process $y(t_i)=H(t_i)x(t_i)+\varepsilon_i$ could easily be considered, without changing formalism, with the same operator $H$ and independent Gaussian noises $\varepsilon_i\sim \cN(0,R_i)$.}
	\begin{equation}
		dy(t)=(H(t)x(t)+h(t))dt+db(t), \quad y(t_0)=y_0\label{eq:2-2}
	\end{equation}
	in which $H(\cdot)\in L^2([t_0,T],\mathcal{L}(\R^{n};\R^{m}))$, $h(\cdot) \in L^1([t_0,T],\R^{m})$. At any final time $T$, the available information is the trajectory $(y(t))_{t\in[t_0,T]}$ and we want to estimate the value of $x(s)$ for $s\in[t_0,T]$. When $s<T$, the best estimate is called the RTS smoother, and when $s=T$, it is the Kalman filter. 
	
	\subsection{BEST ESTIMATE }\label{sec:best-estimate}
	
	From \eqref{eq:2-1} and \eqref{eq:2-2}, we deduce that the processes $x(t)$ and $y(t)$ are Gaussian, with mean $\bar{x}(t)$
	and $\bar{y}(t)$ given by the equations 
	\begin{align}
		\dfrac{d\bar{x}}{dt}(t)&=F(t)\bar{x}(t)+f(t), \quad &\bar{x}(t_0)=x_{0},\label{eq:2-3}\\
		\dfrac{d\bar{y}}{dt}(t)&=H(t)\bar{x}(t)+h(t), \quad &\bar{y}(t_0)=y_0.\label{eq:2-4}
	\end{align}
	We introduce the processes $\widetilde{x}(t)$ and $\widetilde{y}(t)$
	with zero mean, $\widetilde{x}(t):=x(t)-\bar{x}(t)$, $\widetilde{y}(t):=y(t)-\bar{y}(t)$.
	They are solutions of the SDEs
	\begin{align}
		d\widetilde{x}(t)&=F(t)\widetilde{x}(t)dt+G(t)dw(t), \quad &\widetilde{x}(t_0)=\xi,\label{eq:2-5}\\
		d\widetilde{y}(t)&=H(t)\widetilde{x}(t)dt+db(t), \quad &\widetilde{y}(t_0)=0.\label{eq:2-6}
	\end{align}
	The formal problem is to estimate $x(s)$ with the $\sigma$-algebra
	$\mathcal{Y}^{T}=\sigma(y(\tau),0\leq\tau\leq T)$. For the minimum mean square estimator, the solution is
	well-known \citep{Neveu1968ProcessusAG}, it is the conditional expectation:
	\begin{equation}
		\hat{x}(s | T)=E[x(s)|\mathcal{Y}^{T}].\label{eq:2-7}
	\end{equation}
	The random variable $\hat{x}(s | T)$, for $s<T$, is the RTS smoother estimate, while $\hat{x}(T|T)$ is the Kalman filter estimate. However the expression (\ref{eq:2-7}) is not operational and we need more implementable formulas. The most important element we want to exploit is that, thanks to the fact that the processes $x(\cdot)$ and $y(\cdot)$ are Gaussian, the conditional expectation coincides with the best unbiased linear estimate, a.k.a.\ the minimum variance linear estimator. A linear unbiased estimate is characterized by an operator $S(t)\in\mathcal{L}(\R^{m};\R^{n})$. We estimate $x(s)$ by the minimum variance linear estimator $x_{S}(s|T)$ defined by the formula, related to Wiener filtering,
	\begin{equation}
		x_{S}(s | T):=\bar{x}(s)+\int_{t_0}^{T}S_s(t | T)d\widetilde{y}(t)\label{eq:2-8}
	\end{equation}
	which is obviously unbiased. The estimation error $\epsilon_{S}(s | T)$ thus satisfies
	\begin{equation}
		\epsilon_{S}(s | T):=x(s)-x_{S}(s|T)=\widetilde{x}(s)-\int_{t_0}^{T}S_s(t | T)d\widetilde{y}(t).\label{eq:2-9}
	\end{equation}
	The objective is to find $\hat S_s(\cdot | T)$ minimizing the covariance matrix of the error 
	\begin{equation}
		\hat{S}_s(\cdot | T) \in \argmin_{S(\cdot | T)} \Gamma_{S}(s | T)=\bE[\epsilon_{S}(s | T)(\epsilon_{S}(s | T))^{*}].\label{eq:2-10}
	\end{equation}
	This minimization must be interpreted in the sense of the operator norm of positive matrices. The fact that this problem has a solution is a fundamental result of Kalman smoothing and filtering theory. In this work we obtain a new expression for the optimal operator $\hat{S}_s(\cdot | T)$ by  finding a closed-form formula for the following proper covariance function,
	\begin{equation}
		K(s,t | T)=\bE[\epsilon_{\hat S_s}(s | T)(\epsilon_{\hat S_t}(t | T))^{*}]\in\cL(\R^{n,*},\R^n)\label{eq:kernel_def}
	\end{equation}
	We show below in Corollary~\ref{prop3-1} that  $\hat{S}_s(t|T)=K(s,t|T)H^{*}(t)R^{-1}(t)$. From a Gaussian process perspective, $K$ can be interpreted as the posterior covariance given $\widetilde{y}$ with the prior \eqref{eq:2-100} over $x(\cdot)$. The fact that covariances $K$ are also reproducing kernels is well-known since \citet{Parzen_1961}. \citet[Chapter 2]{berlinet04reproducing} summarized the deep relations between positive semidefinite kernels and stochastic processes, albeit only in the real-valued case. \citet{kanagawa2018} also reviewed the connections between kernel methods and Gaussian processes. We meet here a similar relation to the classical ``duality'' between estimation and control in the Linear-Quadratic-Gaussian setting, which we will come back to in \Cref{sec:detour_estimator,sec:dual_problem_rkhs}. The kernels have a ``dual'' nature, both deterministic, as Green functions of differential equations, and stochastic, as covariances of second order processes.
	
	\subsection{A DETOUR THROUGH ESTIMATORS}\label{sec:detour_estimator}
	
	In this section only, we consider infinite-dimensional operators to give the high-level idea before introducing the various quantities studied in the next sections. We will proceed formally, following the introduction of \citet[Chapter 4.5]{Bensoussan2018} on finite-dimensional best estimators and the short presentation of \citet[Chapter 2.4]{berlinet04reproducing} of Hilbert spaces generated by a process with finite second order moments. 
	 For such a process $X=(X_t)_{t\in[t_0,T]}$, we write $\cC_X$ its covariance operator. For simplicity, we consider only zero-mean processes $X,Y$ and only quantities having finite second order moments. Setting $\T=[t_0,T]$, the (Bayesian) minimum mean square estimator (MMSE) is defined as
	\begin{equation}\label{eq:mmse}\tag{MMSE}
		\min_{\hat{X}\in L^2(\Omega\times \T,\R^n),\, \Phi \text{ meas.},\, \hat{X}=\Phi(Y)} \bE((X-\hat X)^\top(X-\hat X)).
	\end{equation}
	where $\Phi$ is a measurable function from $L^2(\Omega\times \T,\R^m)$ to $L^2(\Omega\times \T,\R^n)$, whence $\hat{X}$ belongs to the space $\bar \cN(Y)$ of nonlinear functionals  of $Y$. The linear MMSE, which coincides with the minimum variance linear estimator (MVLE), instead restricts the search space to the space $\bar \cL(Y)$ of linear functionals of the process $Y$
	\begin{equation}\label{eq:mvle}\tag{MVLE}
		\min_{\hat{X}\in L^2(\Omega\times \T,\R^n),\, \cS\in \cL(L^2(\Omega\times \T,\R^m),L^2(\Omega\times \T,\R^n)) ,\, \hat{X}=\cS Y} \bE((X-\hat X)^\top(X-\hat X)).
	\end{equation}
	Gaussian process regression for real-valued $x$, i.e.\ $n=1$, further restricts the search space, by introducing the canonical congruence $\psi_Y$  between the process $Y$ and its RKHS $\cH_Y$, i.e.\ the linear isomorphism satisfying $\psi_Y(v^\top Y_t)(s)=\bE[Y_s Y_t^\top]v$ for all $v\in\R^m$,  
	\begin{equation}\label{eq:gpr}\tag{GP-reg}
		\min_{\hat{X}\in L^2(\Omega\times \T,\R^n),\, g\in \cH_Y,\, \hat{X}=\psi_Y^{-1}(g)} \bE((X-\hat X)^\top(X-\hat X)).
	\end{equation}
	By convex duality, introducing a process $\Lambda$ which acts as a Lagrange multiplier, the dual problem of \eqref{eq:mvle} can be written as
	\begin{equation*}
		\max_{\Lambda \in L^2(\Omega\times \T,\R^{n,*})}\min_{\hat{X}\in L^2(\Omega\times \T,\R^n),\, \cS\in \cL(L^2(\Omega\times \T,\R^m),L^2(\Omega\times \T,\R^n))} \bE((X-\hat X)^\top(X-\hat X))+2\scalidx{\Lambda}{\hat{X}-\cS Y}{L^2(\Omega\times \T,\R^n)}.
	\end{equation*}
	Minimizing over $\cS$ imposes $\Lambda \in \bar \cL(Y)^\perp$, the orthogonal space of $\bar \cL(Y)$ in $L^2(\Omega\times \T,\R^n)$. Minimizing over $\hat{X}$ gives $\hat{X}^\top=X^\top-\hat \Lambda$ so
	\begin{equation}\label{eq:mvle-dual2}\tag{MVLE-dual}
		\min_{\Lambda \in \cL(Y)^\perp}\bE((X^\top-\Lambda)^\top(X^\top-\Lambda)).
	\end{equation}
	
	This covector or adjoint process does not exist in the Gaussian process framework because the duality is not taken in $L^2$ but in $\cH_Y$, seen as self-adjoint. For such orthogonal spaces to appear, one has to work in a larger space than $\cH_Y$, such as $L^2$. Similarly, in optimal control, working in the space of absolutely continuous trajectories leads to defining a covector function. Working in the RKHS of linearly controlled trajectories as in \citet{aubin2020hard_control} bypasses the need of such covector.\footnote{Note also that the theory of Gaussian processes (GPs), despite being a subfield of Bayesian statistics, does not require to manipulate a proper likelihood function to solve a regression problem. Uncannily we can draw another parallel with linear-quadratic (LQ) optimal control, where introducing the value function is not necessary to derive the solution. In a nutshell, the GP (resp.\ kernel) approach to estimation (resp.\ control) in the Gaussian (resp.\ LQ) case does not involve a larger space of functions to work in, but restricts the analysis to the space generated by the process (resp.\ RKHS). In our context, we will see that one can formally move from stochastic Bayesian GPs to deterministic frequentist kernels by replacing the Brownian noise $dw(t)$ by a control $u(t)dt$. This  formal change is related to moving from Itô to Stratonovitch calculus and allows to derive exact properties of the stochastic system such as set invariance based on the deterministic counterpart \citep{DaPrato_Frankowska_2004}, introducing a Stratonovitch drift in case of nonlinear diffusion terms.}
	
	Coming back to \eqref{eq:mvle}, the first order optimality condition writes as the orthogonality criterion $\bE((X-\hat X)Y^\top)=0$ which gives that $\hat{\cS}=\cC_{XY}\cC_{Y}^{-1}$ if $\cC_{Y}$ is invertible. If $Y=\cH X +Z$, with invertible $\cC_{Z}=\cR$, then $\cC_{Y}$ has an inverse and, setting $\cC_\epsilon=\cC_X-\cC_{\hat{X},\hat{S}}$, further assuming $\cC_X$ to be invertible and using the Woodbury identity, we obtain that $\hat{\cS}=\cC_\epsilon \cH^\top \cR^{-1}$. We note also that for any deterministic $\bar \lambda(\cdot)\in L^2([t_0,T],\R^{n,* })$, $\hat{\cS}$ also minimizes
	\begin{equation*}
		\min_{\hat{X}\in L^2(\Omega\times \T,\R^n),\, \cS\in \cL(L^2( \T,\R^m),L^2(\T,\R^n)),\, \hat{X}=\cS Y} \bE(\|\bar \lambda(\cdot)^\top(X-\hat X)\|_{L^2}^2).
	\end{equation*}
	Defining the deterministic $v(\cdot)=\cS^\top \bar \lambda(\cdot)$, $\hat{\cS}$ minimizes a ``control'' problem over $v$, with a specific rewriting for $Y=\cH X +Z$,    
	\begin{equation}\label{eq:mvle-3}\tag{MVLE-det}
		\min_{\cS\in \cL(L^2(\T,\R^m),L^2(\Omega\times \T,\R^n)),\, v(\cdot)=\cS^\top \bar \lambda(\cdot)}  \underbrace{\scalidx{\bar \lambda(\cdot)}{\cC_X \bar \lambda(\cdot)}{L^2} + \scalidx{v(\cdot)}{\cC_Y    v(\cdot)}{L^2}-2\scalidx{\bar \lambda(\cdot)}{\cC_{XY} v(\cdot)}{L^2}}_{\scalidx{\bar \lambda(\cdot)-\cH^\top v(\cdot)}{\cC_{X} (\bar \lambda(\cdot)-\cH^\top v(\cdot))}{L^2}+\scalidx{v(\cdot)}{\cR v(\cdot)}{L^2}}.
	\end{equation}
	Another approach would have been to consider a Bayesian posterior estimate when considering jointly Gaussian processes $X,Y$. However, as recalled in \citet[][p.9 and 16]{kanagawa2018}, Bayes' rule is more involved in infinite dimensions and the likelihood may be degenerate as, for $\cC_X$ non-invertible, $X$ does not have a density w.r.t.\ the Lebesgue measure. However, for invertible joint covariance $\cC_{X,Y}$, we have that $\cC_Y$ and $\cC_\epsilon=\cC_X-\cC_{XY}\cC_Y^{-1}\cC_{YX}$ are invertible. Nevertheless, for continuous-time observations, the absence of a Lebesgue measure requires extensive care \citep{DaPrato2006}, leading to a less common viewpoint on reproducing kernels in relation with Gaussian measures and Cameron-Martin spaces \citep[Definition 2.3.4]{lunardi2016}. We bypass these difficulties by expressing the maximum log-likelihood estimator for a realization $y(\cdot)$ of $Y$ as its resulting least squares problem\footnote{In continuous time, $\tilde{y}(t)\in L^2(\T,\R^m)$ ``is reminiscent of the observation process, in fact rather the	derivative of the observation process (which, as we know, does not exist)'' \Citep[p180]{Bensoussan2018}. It is as if we claim to observe the derivative, but $\tilde{y}(t)$ will always appear within integrals.\label{foo:observation_continue}}
	\begin{multline}\label{eq:mvle-ls}\tag{LSE}
		\min_{\hat{x}(\cdot)\in L^2(\T,\R^n)} \scalidx{(\hat{x}(\cdot), y(\cdot))}{\cC_{X,Y}^{-1} (\hat{x}(\cdot), y(\cdot))}{L^2(\Omega\times \T,\R^n\times \R^m)}\\
		= \underbrace{\scalidx{\hat{x}(\cdot)-\cC_{XY}\cC_Y^{-1}y(\cdot)}{\cC_{\epsilon}^{-1} (\hat{x}(\cdot)-\cC_{XY}\cC_Y^{-1}y(\cdot))}{L^2}+\scalidx{y(\cdot)}{\cC_Y^{-1}    y(\cdot)}{L^2}}_{\scalidx{\hat{x}(\cdot)}{\cC_X^{-1}    \hat{x}(\cdot)}{L^2}+\scalidx{y(\cdot)-\cH \hat{x}(\cdot)}{\cR^{-1} (y(\cdot)-\cH \hat{x}(\cdot))}{L^2}}.
	\end{multline}
	\begin{remark}[Stochastic and deterministic dual or equivalent problems]
			\citet[Chapter 15, Table 15.1, p568]{Kailath2000} describes the four problems we obtained as either ``dual'' or ``equivalent''. Here we intend duality strictly in the Fenchel context, obtained by permutation of max and min. The problems \eqref{eq:mvle}-\eqref{eq:mvle-dual2} are indeed stochastic dual problems. That \eqref{eq:mvle-ls}-\eqref{eq:mvle-3} are deterministic dual problems as claimed in \citet[p59]{Bensoussan2018} is not so straightforward (this is justified through dual bases in \citet{Kailath2000}). We will return to this in \Cref{sec:dual_problem_rkhs} and prove that when expressed on RKHSs, there is indeed a Fenchel duality. Similarly, that \eqref{eq:mvle}-\eqref{eq:mvle-ls} and \eqref{eq:mvle-dual2}-\eqref{eq:mvle-3} are ``equivalent'' stochastic and deterministic problems, as defined in \citet[Section 15.3]{Kailath2000}, will get clearer from the fact that both problems share the same reproducing kernel (see \Cref{sec:rkhs}). Note that \eqref{eq:mvle-ls}-\eqref{eq:mvle-3} do not have the same assumptions concerning the invertibility of $\cC_{X,Y}$ and of $\cC_{Y}$, a limitation which we will not meet when using kernels.
	\end{remark}
	The problem with the approaches developed so far is that manipulating abstract covariances and continuous linear maps over $L^2$ can be cumbersome and is not operable for non-discrete time. Instead, Laurent Schwartz's kernel theorem allows us to consider all these operators as kernel integral operators. Thus, in the following, we will see formulas very similar to the ones derived in this section, but written explicitly on the time axis and with recursive versions in the spirit of Kalman filtering. 
	Formally one can consider in \eqref{eq:mvle-3} the test function $\bar \lambda(\cdot)=\delta_s(\cdot) \bar \lambda$ for some $\bar \lambda\in \R^{n,*}$ to recover \eqref{eq:2-10}, but we will see that choosing $\bar \lambda(\cdot)$ based on the adjoint process $\hat \Lambda$ will allow us to avoid having to manipulate $\cC_X$ in \eqref{eq:mvle-3}.
	
	\subsection{A RELATED OPTIMAL CONTROL PROBLEM}
	We return to the formulas of \Cref{sec:best-estimate} and start by expressing $\scalh{\Gamma_{S}(s | T)\bla}{\bla}$ more explicitly,
	\begin{align*}
		\scalh{\bla}{\epsilon_{\hat S_s}(s | T)}&=\scalh{\bla}{\widetilde{x}(s)}-\int_{t_0}^{T}\scalh{S^{*}_s(t | T)\bla}{d\widetilde{y}(t)}=\scalh{\bla}{\widetilde{x}(s)}-\int_{t_0}^{T}\scalh{S^{*}_s(t | T)\bla}{H(t)\widetilde{x}(t)dt+db(t)}.
	\end{align*}
	We seek an expression where $\widetilde{x}$ does not appear. We thus introduce the following adjoint equation over an adjoint variable $\lambda_s$, a.k.a.\ the information vector,
	\begin{equation}
		-\dfrac{d\lambda_s}{dt}=F^{*}(t)\lambda_s(t)-H^{*}(t)S^{*}_s(t | T)\bla, \quad \lambda_s(T)=\left|\begin{array}{cc}
			0 & \text{if}\;s<T\\
			\bla & \text{if}\;s=T
		\end{array}\right., \quad \lambda_{s}(s)-\lambda_{s}(s^+)=\bla,\text{if}\:s<T,\label{eq:2-11}
	\end{equation}
	where $\lambda_{s}(s^+)=\lim_{h\rightarrow s^+}\lambda_{s}(h)$ with the convention $\lambda_{s}(T^+)=0$. Then a simple calculation shows that
	\begin{align*}
		\scalh{\bla}{\epsilon_{\hat S_s}(s | T)}&=\scalh{\bla}{\widetilde{x}(s)}+\int_{t_0}^{T}\scalh{-\dfrac{d\lambda}{dt}-F^{*}(t)\lambda_s(t)}{\widetilde{x}(t)dt}-\int_{t_0}^{T}\scalh{S^{*}_s(t | T)\bla}{db(t)}\\
		&\hspace{-2cm}=\scalh{\lambda_{s}(s)-\lambda_{s}(s^+)}{\widetilde{x}(s)}+\int_{t_0}^{T}\scalh{-\dfrac{d\lambda}{dt}}{\widetilde{x}(t)dt}+\int_{t_0}^{T}\scalh{\lambda_s(t)}{-d\widetilde{x}(t)+G(t)dw(t)}-\int_{t_0}^{T}\scalh{S^{*}_s(t | T)\bla}{db(t)}\\
		&\hspace{-2cm}=\underbrace{\scalh{\lambda_{s}(s)-\lambda_{s}(s^+)}{\widetilde{x}(s)}-[\scalh{\lambda_{s}(\tau)}{\widetilde{x}(\tau)}]_{t_0}^s-[\scalh{\lambda_{s}(\tau)}{\widetilde{x}(\tau)}]_{s^+}^T}_{=\scalh{\lambda_{s}(t_0)}{\xi}}+\int_{t_0}^{T}\scalh{G^*(t)\lambda_s(t)}{dw(t)}-\int_{t_0}^{T}\scalh{S^{*}_s(t | T)\bla}{db(t)}
	\end{align*}
	and therefore 
	\begin{equation}
		\scalh{\Gamma_{S}(s|T)\bla}{\bla}=\scalh{\Pi_{0}\lambda_s(t_0)}{\lambda_s(t_0)}+\int_{t_0}^{T}\scalh{G(t)Q(t)G^{*}(t)\lambda_s(t)}{\lambda_s(t)}dt+\int_{t_0}^{T}\scalh{R(t)S^{*}_s(t|T)\bla}{S^{*}_s(t|T)\bla}dt\label{eq:2-12}
	\end{equation}
	More generally, beyond linear feedbacks ($S^{*}_s(\cdot | T)\bla$), for a general control input $v(\cdot)$ with $R(\cdot)^{\frac{1}{2}}v(\cdot)\in L^2(t_0,T;\R^m)$, it is natural to consider the following control problem extending \eqref{eq:2-11} and corresponding to \eqref{eq:mvle-3}
	\begin{gather}
		-\dfrac{d\lambda_s}{dt}=F^{*}(t)\lambda_s(t)+H^{*}(t)v(t), \, \lambda_s(T)=\left|\begin{array}{cc}
			0 & \text{if}\;s<T\\
			\bla & \text{if}\;s=T
		\end{array}\right., \quad
		\lambda_{s}(s)-\lambda_{s}(s^+)=\bla,\text{if}\:s<T;\label{eq:2-13}\\
		J(v(\cdot))=\scalh{\Pi_{0}\lambda_s(t_0)}{\lambda_s(t_0)}+\int_{t_0}^{T}\scalh{G(t)Q(t)G^{*}(t)\lambda_s(t)}{\lambda_s(t)}dt+\int_{t_0}^{T}\scalh{R(t)v(t)}{v(t)}dt.\label{eq:2-14}
	\end{gather}
	The classical way to solve such a problem is through Pontryagin's Maximum Principle \citep[see e.g.][Chapter 10]{Bensoussan2018}, i.e.\ by using Lagrange-Fenchel duality, which effectively leads to a Hamiltonian system, taking the form of a two-point boundary value problem:
	\begin{align}
		\dfrac{d\hat{\gamma}_s}{dt}&=F(t)\hat{\gamma}_s(t)-G(t)Q(t)G^{*}(t)\hat{\lambda}_s(t)\label{eq:2-15}\\
		-\dfrac{d\hat{\lambda}_s}{dt}&=F^{*}(t)\hat{\lambda}_s(t)+H^{*}(t)R^{-1}(t)H(t)\hat{\gamma}_s(t) \nonumber \\
		\hat{\gamma}_s(t_0)&=-\Pi_{0}\hat{\lambda}_s(t_0), \quad \lambda_s(T)=\left|\begin{array}{cc}
			0 & \text{if}\;s<T\\
			\bla & \text{if}\;s=T
		\end{array}\right., \quad
		\lambda_{s}(s)-\lambda_{s}(s^+)=\bla,\text{if}\:s<T. \nonumber
	\end{align}
	and where the optimal control of problem \eqref{eq:2-13}-\eqref{eq:2-14} is given by
	\begin{equation}
		\hat{v}_s(t)=R^{-1}(t)H(t)\hat{\gamma}_s(t).\label{eq:2-16}
	\end{equation}
	Since the pair $(\hat{\lambda}_s(t),\hat{\gamma}_s(t))$ depends linearly on $\bla$, we obtain immediately that there exists a single $\widehat{S}_s(\cdot | T)$ which minimizes (\ref{eq:2-12}) for any $\bla$, namely 
	\begin{equation}
		\widehat{S}_s^{*}(t | T)\bla=-R^{-1}(t)H(t)\hat{\gamma}_s(t).\label{eq:2-17}
	\end{equation}
	This suggests that one has to solve \eqref{eq:2-15} to get an expression of $\widehat{S}^{*}(s,t | T)$.
	
	\subsection{A GENERAL TWO-POINT BOUNDARY VALUE PROBLEM}\label{sec:twoPtsBdry}
	At this stage, rather than coping with the difficulty of dealing with jumps as in \eqref{eq:2-15}, we consider another two-point boundary value problem similar to (\ref{eq:2-15}) but which crucially does not depend on $s$. We also had not so far considered the possibility of having a Gaussian prior with covariance $\Sigma_T\in \cL(\R^n,\R^{n,*})$ on the terminal information $\lambda_s(T)$, so we introduce it here for greater generality. This term will act a special weight on the terminal point $x(T)$ in relation with a terminal cost in optimal control. Consider the pair $(\widehat{\mu}(t),\widehat{\nu}(t))$ solution of the coupled system
	\begin{align}
		\dfrac{d\hat{\mu}}{dt}&=F(t)\hat{\mu}(t)-G(t)Q(t)G^{*}(t)\hat{\nu}(t)+l_\mu(t)\label{eq:twoPtsBdry}\\
		-\dfrac{d\hat{\nu}}{dt}&=F^{*}(t)\hat{\nu}(t)+H^{*}(t)R^{-1}(t)H(t)\hat{\mu}(t)-l_\nu(t) \nonumber \\
		\hat{\mu}(t_0)&=-\Pi_{0}\hat{\nu}(t_0),\;\hat{\nu}(T)=\Sigma_T \hat{\mu}(T). \nonumber
	\end{align}
	where $\hat{\mu}$ (resp.\ $\hat{\nu}$) plays the role of $\gamma_s$ (resp.\ $\lambda_s$) and $l_\mu(\cdot)\in L^{2}(t_0,T;\R^{n})$ and $l_\nu(\cdot)\in L^{2}(t_0,T;\R^{n,*})$ are two test functions acting as perturbations of the differential equations. Notice that, for  $l_\mu(\cdot)\equiv 0$ and $l_{\nu,s}(\tau)=-\bla \delta_s(\tau)$, the system \eqref{eq:2-15} with a jump condition corresponds precisely to the two-point boundary system \eqref{eq:twoPtsBdry}. The important result is the following
	\begin{prop}\label{prop2-1} For $\Sigma_T=0$, $l_\mu(\cdot)\equiv 0$ and $l_\nu(t)=H^{*}(t)R^{-1}(t)g(t)$ with $g(\cdot)\in L^{2}(t_0,T;\R^{m})$, we have the formula 
		\begin{equation}
			\hat{\mu}(s)=\int_{t_0}^{T}\hat{S}_s(t | T)g(t)dt.\label{eq:2-200}
		\end{equation}
	\end{prop}
	
	\begin{proof}
		From (\ref{eq:2-17}) we have 
		\[
		\scalh{\bla}{\int_{t_0}^{T}\widehat{S}_s(t | T)g(t)dt}=-\int_{t_0}^{T}\scalh{\hat{\gamma}_s(t)}{H^{*}(t)R^{-1}(t)g(t)}\,dt.
		\]
		We use the second equation of (\ref{eq:twoPtsBdry}) and proceed with integration
		by parts between the systems (\ref{eq:2-15}) and (\ref{eq:twoPtsBdry}). Note that $\hat{\mu}$, $\hat{\nu}$, $\hat{\gamma}_s$ are all continuous functions, but that $\hat{\lambda}_s$ is not, so that we have to split some integrals to take care of the singularity at $s$ as we did before to derive \eqref{eq:2-12}.
		\begin{align*}
			&\scalh{\bla}{\int_{t_0}^{T}\widehat{S}_s(t | T)g(t)dt}=\int_{t_0}^{T}\scalh{\hat{\gamma}_s(t)}{-\dfrac{d\hat{\nu}}{dt}-F^{*}(t)\hat{\nu}(t)-H^{*}(t)R^{-1}(t)H(t)\hat{\mu}(t)}dt\\
			&=\int_{t_0}^{T}\scalh{\hat{\gamma}_s(t)}{-\dfrac{d\hat{\nu}}{dt}}dt-\int_{t_0}^{T}\scalh{\dfrac{d\hat{\gamma}_s}{dt}+G(t)Q(t)G^{*}(t)\hat{\lambda}_s(t)}{\hat{\nu}(t)}dt+\int_{t_0}^{T}\scalh{\dfrac{d\hat{\lambda}_s}{dt}+F^{*}(t)\hat{\lambda}_s(t)}{\hat{\mu}(t)}dt\\
			&=-[\scalh{\hat{\gamma}_s(t)}{\hat{\nu}(t)}]_{t_0}^T+\int_{t_0}^{T}\scalh{\hat{\lambda}_s(t)}{\dfrac{d\hat{\mu}}{dt}-F(t)\hat{\mu}(t)}dt+\int_{t_0}^{T}\scalh{\dfrac{d\hat{\lambda}_s}{dt}+F^{*}(t)\hat{\lambda}_s(t)}{\hat{\mu}(t)}dt\\
			&=-\scalh{\hat{\gamma}_s(T)}{\underbrace{\hat{\nu}(T)}_{=\Sigma_{T}\hat{\mu}(T)=0}}+\scalh{\underbrace{\hat{\gamma}_s(t_0)}_{=-\Pi_{0}\hat{\lambda}_s(t_0)}}{\hat{\nu}(t_0)}+\scalh{\underbrace{\hat{\lambda}_s(T)}_{=0}}{\hat{\mu}(T)}\1_{s<T}-\scalh{\hat{\lambda}_s(s^+)-\hat{\lambda}_s(s)}{\hat{\mu}(s)}-\scalh{\hat{\lambda}_s(t_0)}{\underbrace{\hat{\mu}(t_0)}_{=-\Pi_{0}\hat{\nu}(t_0)}}\\
			&=\scalh{\hat{\lambda}_s(s)-\hat{\lambda}_s(s^+)}{\hat{\mu}(s)}=\scalh{\bla}{\hat{\mu}(s)}.
		\end{align*}
	\end{proof}
	
	To derive another formula for $\widehat{S}_s(t | T)$, we are now going to find another matrix satisfying \eqref{eq:2-200} for all $g(\cdot)$ with $\hat{\mu}(\cdot)$ defined as per \eqref{eq:twoPtsBdry}.
	
	\section{SOLUTION THROUGH RICCATI EQUATIONS AND KERNELS}\label{sec:riccati-kernels}
	
	\subsection{DERIVING THE KERNELS}
	When $l_\mu(\cdot)\equiv 0$ and $l_\nu(\cdot)\equiv 0$, the classical approach to solve \eqref{eq:twoPtsBdry} is by variation of constants, introducing two matrices $\Sigma(t)$ and $\Pi(t)$ satisfying $\hat{\mu}(t)=-\Pi(t)\hat{\nu}(t)$ and $\hat{\nu}(t)=\Sigma(t) \hat{\mu}(t)$. It is then straightforward to show that they must satisfy two (dual) differential Riccati equations
	\begin{align}
		-\dfrac{d}{dt}\Sigma&=\Sigma(t)F(t)+F^{*}(t)\Sigma(t)-\Sigma(t)G(t)Q(t)G^{*}(t)\Sigma(t)+H^{*}(t)R^{-1}(t)H(t), \quad &\Sigma(T)=\Sigma_T;\label{eq:2-20}\\
		\dfrac{d}{dt}\Pi&=F(t)\Pi(t)+\Pi(t)F^{*}(t)-\Pi(t)H^{*}(t)R^{-1}(t)H(t)\Pi(t)+G(t)Q(t)G^{*}(t),\quad &\Pi(t_0)=\Pi_{0}.\label{eq:5-110}
	\end{align}
	The information filter matrix $\Sigma(t)$ satisfies a backward equation depending on $T$ whereas the estimation filter matrix $\Pi(t)$ satisfies a forward one and depends on $t_0$. The solutions for these equations exist at all times on $[t_0,T]$. This is a classical result whose proof we give in Lemma~\ref{lem3-1} in Appendix A.1.\footnote{If $\Pi_{0}$ is invertible, then so is $\Pi(t)$. Indeed the Riccati equation preserves positive definiteness \citep[see e.g.][Example 16.3.4, p629]{Kailath2000}. Then $\Pi(t)^{-1}$ also satisfies \eqref{eq:2-20}, and, if we choose $\Sigma_T=\Pi(T)^{-1}$, then we can identify the matrices $\Sigma(t)$ and $\Pi(t)^{-1}$.}\\

	When $l_\mu(\cdot)\not\equiv 0$ or $l_\nu(\cdot)\not\equiv 0$, we are going instead to follow a Green kernel approach, related to Wiener filtering, and look for kernel integral operators satisfying respectively:\footnote{An alternative scheme would be to consider exponentials of the Hamiltonian matrix as in \citet[Chapter 5]{speyer2010primer}. We use and recall this approach in the Appendix A.2 to give a numerical method to evaluate the kernel $K(s,t|T)$.}
	\begin{align}
		\hat{\mu}(s)=\int_{t_0}^{T}K(s,t | T)l_\nu(t)dt \quad &\text{ for } l_\mu(\cdot)\equiv 0, \label{eq:Wiener_K}\\
		\hat{\nu}(s)=\int_{t_0}^{T}\Lambda(s,t | T)l_\mu(t)dt  \quad &\text{ for } l_\nu(\cdot)\equiv 0. \label{eq:Wiener_L}
	\end{align}
	These formulas may recall the ones appearing in the innovations approach \citep[Section 16.4.2]{Kailath2000}; however we do not consider cross-covariances between trajectories and observations. We will show (see Proposition~\ref{prop:covar_BF_formula} below) that $K$ satisfying \eqref{eq:Wiener_K} is precisely the covariance of the optimal error \eqref{eq:kernel_def}. On the other hand $\Lambda$ as in \eqref{eq:Wiener_L} is related to observability problems. To solve \eqref{eq:twoPtsBdry} in general, we can compute the distorsion w.r.t.\ the  solutions with null perturbations, introducing two variables:
	\begin{equation}
		r(t)=\hat{\mu}(t)+\Pi(t)\hat{\nu}(t),\; \eta(t)=\hat{\nu}(t)-\Sigma(t) \hat{\mu}(t).\label{eq:error_edo}
	\end{equation}
	From \eqref{eq:twoPtsBdry}, \eqref{eq:2-20} and \eqref{eq:5-110}, we deduce that $r(\cdot)$ and $\eta(\cdot)$ satisfy the following differential equations\footnote{We give below the computation for \eqref{eq:edo_r}, eq.\eqref{eq:edo_eta} is obtained similarly \begin{equation*}
			dr/dt=F\hat{\mu}-GQG^{*}\hat{\nu}+l_\mu+(F\Pi+\Pi F^{*}-\Pi H^{*}R^{-1}H\Pi+GQG^{*})\hat{\nu}-\Pi(F^{*}\hat{\nu}+H^{*}R^{-1}H\hat{\mu}-l_\nu).
	\end{equation*}}
	\begin{align}
		\dfrac{d}{dt}r&=(F(t)-\Pi(t)H^{*}(t)R^{-1}(t)H(t))r(t)+\Pi(t)l_\nu(t)+l_\mu(t),\quad &r(t_0)=0;\label{eq:edo_r}\\
		-\dfrac{d}{dt}\eta&=(F^*(t)-\Sigma(t)G(t)Q(t)G^{*}(t))\eta(t)+\Sigma(t)l_\mu(t)-l_\nu(t),\quad &\eta(T)=0.\label{eq:edo_eta}
	\end{align}
	This suggests to introduce semigroups associated with the matrix function $F(t)-G(t)Q(t)G^{*}(t)\Sigma(t)$ denoted $\Phi_{F,\Sigma}(s,t)$ (resp.\ $F(s)-\Pi(s)H^{*}(s)R^{-1}(s)H(s)$, denoted $\Phi_{F,\Pi}(s,t)$)
	\begin{align*}
	\dfrac{d}{d\tau}\Phi_{F,\Sigma}(\tau,t)&=(F(\tau)-G(\tau)Q(\tau)G^{*}(\tau)\Sigma(\tau))\Phi_{F,\Sigma}(\tau,t),\quad &&\Phi_{F,\Sigma}(t,t)=\Id;\\
	-\dfrac{d}{dt}\Phi_{F,\Sigma}^{*}(\tau,t)&=(F^{*}(t)-\Sigma(t)G(t)Q(t)G^{*}(t))\Phi_{F,\Sigma}^{*}(\tau,t),\quad &&\Phi_{F,\Sigma}(\tau,\tau )=\Id;\\
	\dfrac{d}{d\tau }\Phi_{F,\Pi}(\tau,t)&=(F(\tau )-\Pi(\tau )H^{*}(\tau )R^{-1}(\tau )H(\tau ))\Phi_{F,\Pi}(\tau,t),\quad &&\Phi_{F,\Pi}(t,t)=\Id;\\
	-\dfrac{d}{dt}\Phi_{F,\Pi}^{*}(\tau,t)&=(F^{*}(t)-H^{*}(t)R^{-1}(t)H(t)\Pi(t))\Phi_{F,\Pi}^{*}(\tau,t),\quad &&\Phi_{F,\Pi}^{*}(\tau,\tau )=\Id.
\end{align*}
	Note that if $Q(\cdot)\equiv 0$ and $\Pi_0=0$, then $\Pi(\cdot)\equiv 0$ and $\Phi_{F,\Pi}(t,s)=\Phi_{F}(t,s)$, the semi-group associated with the operator $F(t)$, i.e.\ $\partial_t \Phi_{F}(t,s)=F(t)\Phi_{F}(t,s)$ and $\Phi_{F}(t,t)=\Id$. Similarly, if $H(\cdot)\equiv 0$ and $\Sigma_T=0$, then $\Phi_{F,\Sigma}(t,s)=\Phi_{F}(t,s)$. In the general case, we have
	
	\begin{thm}\label{thm:main_noyaux}
		The kernels $K$ and $\Lambda$ satisfying \eqref{eq:Wiener_K} and \eqref{eq:Wiener_L} are given by
		\begin{align}
			K(s,t | T)&=\Phi_{F,\Sigma}(s,t_0)\Pi_{0}^{\frac{1}{2}}(\Id+\Pi_{0}^{\frac{1}{2}}\Sigma(t_0)\Pi_{0}^{\frac{1}{2}})^{-1}\Pi_{0}^{\frac{1}{2}}\Phi_{F,\Sigma}^{*}(t,t_0)+\int_{t_0}^{\min(s,t)}\Phi_{F,\Sigma}(s,\tau)G(\tau)Q(\tau)G^{*}(\tau)\Phi_{F,\Sigma}^{*}(t,\tau)d\tau\label{eq:3-3};\\
			\Lambda(s,t|T)&=\Phi_{F,\Pi}^{*}(T,s)\Sigma_T^{\frac{1}{2}}(\Id+\Sigma_T^{\frac{1}{2}}\Pi(T)\Sigma_T^{\frac{1}{2}})^{-1}\Sigma_T^{\frac{1}{2}}\Phi_{F,\Pi}(T,t)+\int_{\max(s,t)}^{T}\Phi_{F,\Pi}^{*}(\tau,s)H^{*}(\tau)R^{-1}(\tau)H(\tau)\Phi_{F,\Pi}(\tau,t)d\tau\label{eq:5-16}.
		\end{align}
	\end{thm}
	\tb{Remark.} The matrix inverses we consider always exist since $\Pi_{0}$, $\Sigma_T$, $\Sigma(t_0)$ and $\Pi(T)$ are all positive semidefinite. Note that if $\Pi_{0}$ is invertible, then the first term in $K$ boils down to $\Phi_{F,\Sigma}(s,t_0)(\Pi_{0}^{-1}+\Sigma(t_0))^{-1}\Phi_{F,\Sigma}^{*}(t,t_0)$. If $\Sigma_T=0$, then the first term in $\Lambda$ vanishes. If furthermore $H(\cdot)\equiv 0$, the observations are then pure Brownian noise independent from the state, and thus dispensable. We hence recover the formula of \citet[Section 6.4, p.89]{Srkk2019} for the (prior) covariance of a linear SDE. In \citet[Example 3.3]{Kailath2000}, the case $G(\cdot)\equiv 0$ was further discussed (see also \Cref{sec:dualRKHS}) as it makes the second term in $K$ vanish. Both kernels satisfy a Hermitian symmetry, i.e.\ $K(s,t|T)=K(t,s|T)^*$. In Appendix A.2, we give a numerical method to evaluate the kernel $K(s,t|T)$ through the Hamiltonian semigroup. For time-invariant systems, this boils down to matrix exponentials.
	\begin{proof}
		Let us first identify $K$. By the variation of constants formula, \eqref{eq:edo_eta} yields
		\begin{equation}
			\hat{\nu}(t)-\Sigma(t) \hat{\mu}(t)=\eta(t)=\int_{t}^T \Phi^*_{F,\Sigma}(\tau,t)(\Sigma(\tau)l_\mu(\tau)-l_\nu(\tau)) d\tau.\label{eq:eta_explicite}
		\end{equation}
		Consequently \eqref{eq:twoPtsBdry} becomes
		\begin{equation}
			\dfrac{d\hat{\mu}}{dt}=F(t)\hat{\mu}(t)-G(t)Q(t)G^{*}(t)\left(\Sigma(t) \hat{\mu}(t)+\int_{t}^T \Phi^*_{F,\Sigma}(\tau,t)(\Sigma(\tau)l_\mu(\tau)-l_\nu(\tau)\right) d\tau)+l_\mu(t).
		\end{equation}
		Now consider the case where $l_\mu(\cdot)\equiv 0$. Introducing the square root $\Pi_{0}^{\frac{1}{2}}$, the variation of constants formula gives
		\begin{align}
			\hat{\mu}(s)&=\Phi_{F,\Sigma}(s,t_0)\hat{\mu}(t_0)+\int_{t_{0}}^s \Phi_{F,\Sigma}(s,t) G(t)Q(t)G^{*}(t)\int_{t}^T \Phi^*_{F,\Sigma}(\tau,t)l_\nu(\tau)d\tau dt\nonumber \\
			&=-\Phi_{F,\Sigma}(s,t_0)\Pi_{0}^{\frac{1}{2}}\Pi_{0}^{\frac{1}{2}}\hat{\nu}(t_0)+\int_{t_{0}}^T \int_{t_0}^{\min(s,t)}\Phi_{F,\Sigma}(s,\tau)G(\tau)Q(\tau)G^{*}(\tau)\Phi_{F,\Sigma}^{*}(t,\tau)l_\nu(t)d\tau dt.\label{eq:mu_explicite}
		\end{align}
		To compute $\hat{\nu}(t_0)$, we use \eqref{eq:eta_explicite}
		\begin{equation}
			(\Id+\Pi_{0}^{\frac{1}{2}}\Sigma(t_0)\Pi_{0}^{\frac{1}{2}})\Pi_{0}^{\frac{1}{2}}\hat{\nu}(t_0)=\Pi_{0}^{\frac{1}{2}}(\hat{\nu}(t_0)+\Sigma(t_0)\Pi_0 \hat{\nu}(t_0))=\Pi_{0}^{\frac{1}{2}} \eta(t_0)=-\Pi_{0}^{\frac{1}{2}}\int_{0}^T \Phi^*_{F,\Sigma}(\tau,t_0)l_\nu(\tau) d\tau.\label{eq:nu0_explicite}
		\end{equation}
		Inserting \eqref{eq:nu0_explicite} in \eqref{eq:mu_explicite}, we then identify $K$ through \eqref{eq:Wiener_K} which yields \eqref{eq:3-3}. Now for $\Lambda$, the procedure is the same, \eqref{eq:edo_r} and \eqref{eq:twoPtsBdry} yield
		\begin{gather}
			\hat{\mu}(t)+\Pi(t)\hat{\nu}(t)=r(t)=\int_{t_{0}}^t \Phi_{F,\Pi}(t,\tau)(\Pi(\tau)l_\nu(\tau)+l_\mu(\tau)) d\tau, \label{eq:r_explicite}\\
			-\dfrac{d\hat{\nu}}{dt}=F^{*}(t)\hat{\nu}(t)+H^{*}(t)R^{-1}(t)H(t)\left(-\Pi(t)\hat{\nu}(t)+\int_{t_{0}}^t \Phi_{F,\Pi}(t,\tau)(\Pi(\tau)l_\nu(\tau)+l_\mu(\tau)) d\tau\right)-l_\nu(t).\label{eq:edo_nu_explicite}
		\end{gather}
		Now consider the case where $l_\nu(\cdot)\equiv 0$,
		\begin{align}
			\hat{\nu}(s)&=\Phi_{F,\Pi}^{*}(T,s)\hat{\nu}(T)+\int_{t}^T\Phi_{F,\Pi}^{*}(t,s) H^{*}(t)R^{-1}(t)H(t)\int_{t_{0}}^t \Phi_{F,\Pi}(t,\tau)l_\mu(\tau)d\tau dt\nonumber \\
			&=\Phi_{F,\Pi}^{*}(T,s)\Sigma_T^{\frac{1}{2}}\Sigma_T^{\frac{1}{2}}\hat{\mu}(T)+\int_{t_{0}}^T \int_{\max(s,t)}^{T}\Phi_{F,\Pi}^{*}(\tau,s)H^{*}(\tau)R^{-1}(\tau)H(\tau)\Phi_{F,\Pi}(\tau,t)d\tau l_\mu(t) dt.\label{eq:nu_explicite}
		\end{align}
		To compute $\hat{\mu}(T)$, we use \eqref{eq:r_explicite}
		\begin{equation}
			(\Id+\Sigma_T^{\frac{1}{2}}\Pi(T)\Sigma_T^{\frac{1}{2}})\Sigma_T^{\frac{1}{2}}\hat{\mu}(T)=\Sigma_T^{\frac{1}{2}}(\hat{\mu}(T)+\Pi(T)\Pi_0 \hat{\mu}(T))=\Sigma_T^{\frac{1}{2}} r(T)=-\Sigma_T^{\frac{1}{2}}\int_{0}^T \Phi_{F,\Pi}(\tau,t_0)l_\nu(\tau) d\tau.\label{eq:muT_explicite}
		\end{equation}
		Inserting \eqref{eq:muT_explicite} in \eqref{eq:nu_explicite}, we then identify $\Lambda$ through \eqref{eq:Wiener_L} which yields \eqref{eq:5-16}.
		
	\end{proof}
	By Proposition~\ref{prop2-1}, we have shown \eqref{eq:2-200} which matches with \eqref{eq:Wiener_K} and should hold for all test functions $g(\cdot)\in L^{2}(t_0,T;\R^{m})$, so we deduce immediately that: 
	\begin{Corollary}
		\label{prop3-1} If $\Sigma_T=0$, then, for all $s,t\in[t_0,T]$, we have the formula
		\begin{equation}
			\hat{S}_s(t|T)=K(s,t|T)H^{*}(t)R^{-1}(t).\label{eq:3-4}
		\end{equation}
	\end{Corollary}
	
	\subsection{CONTROL INTERPRETATION AS IMPULSE RESPONSE}\label{sec:impulse-response}
	
	A natural quantity to consider in \eqref{eq:Wiener_K} is the impulse-response for $l_\nu(\tau)=z\delta_t(\tau)$ and some $z\in\R^{n,*}$, which gives $x_{zt}(s)=K(s,t | T)z$. Now \eqref{eq:twoPtsBdry} then writes as
	\begin{align}
		\dfrac{d}{d\tau}x_{zt}(\tau)&=F(\tau)x_{zt}(\tau)-G(\tau)Q(\tau)G^{*}(\tau)\nu_{zt}(\tau)\label{eq:3-13} \\
		-\dfrac{d}{d\tau}\nu_{zt}(\tau)&=F^{*}(\tau)\nu_{zt}(\tau)+H^{*}(\tau)R^{-1}(\tau)H(\tau)x_{zt}(\tau)-z\delta_t(\tau) \nonumber \\
		x_{zt}(t_0)&=-\Pi_{0}\nu_{zt}(t_0),\;\nu_{zt}(T)=\Sigma_T x_{zt}(T),\;\nonumber 
	\end{align}
	where the discontinuity occurs in the second line rather than the first. This, we are going to reverse through a change of variables. Considering the function $\chi_{zt}(\tau)$ solution of 
	\begin{equation}
		-\dfrac{d}{d\tau}\chi_{zt}(\tau)=F^{*}(\tau)\chi_{zt}(\tau),\tau<t, \quad \chi_{zt}(t)=z\label{eq:3-14}
	\end{equation}
	we can write formally that $\chi_{zt}(\tau)\1_{\tau<t}$ is solution of
	\begin{equation}
		-\dfrac{d}{d\tau}(\chi_{zt}(\tau)\1_{\tau<t})=F^{*}(\tau)(\chi_{zt}(\tau)\1_{\tau<t})+z\delta_t(\tau), \quad \chi_{zt}(T)\1_{T<t}=0\label{eq:3-15}
	\end{equation}
	With $l_\nu(\tau)=z\delta_t(\tau)$, we define $q_{zt}(\tau)=\nu_{zt}(\tau)+\chi_{zt}(\tau)\1_{\tau<t}$,  then from (\ref{eq:3-13}) the pair $x_{zt}(\tau),q_{zt}(\tau)$ is solution of the system 
	\begin{align}
		\dfrac{d}{d\tau}x_{zt}(\tau)&=F(\tau)x_{zt}(\tau)-G(\tau)Q(\tau)G^{*}(\tau)q_{zt}(\tau)+G(\tau)Q(\tau)G^{*}(\tau)\chi_{zt}(\tau)\1_{\tau<t}\label{eq:3-16} \\
		-\dfrac{d}{d\tau}q_{zt}(\tau)&=F^{*}(\tau)q_{zt}(\tau)+H^{*}(\tau)R^{-1}(\tau)H(\tau)x_{zt}(\tau) \nonumber \\
		x_{zt}(t_0)&=-\Pi_{0}q_{zt}(t_0)+\Pi_{0}\chi_{zt}(t_0),\; q_{zt}(T)=\Sigma_T x_{zt}(T).\nonumber 
	\end{align}
	Similarly to \eqref{eq:2-14}, we interpret the system (\ref{eq:3-16}) as the necessary and sufficient optimality condition of the following control problem 
	\begin{equation}
		\dfrac{d}{d\tau}\zeta(\tau)=F(\tau)\zeta(\tau)+G(\tau)Q^{\frac{1}{2}}(\tau)u(\tau)+G(\tau)Q(\tau)G^{*}(\tau)\chi_{zt}(\tau)\1_{\tau<t}, \quad \zeta(t_0)=\Pi_{0}^{\frac{1}{2}}\xi+\Pi_{0}\chi_{zt}(t_0)\label{eq:3-18}
	\end{equation}
	in which $u(\cdot)$ and $\xi$ are controls, and the objective to minimize is
	\begin{equation}
		J_x(\xi,u(\cdot))=\|\xi\|^{2}+\scalh{\Sigma_T\zeta(T)}{\zeta(T)}+\int_{t_0}^{T}\|u(\tau)\|^{2}d\tau+\int_{t_0}^{T}\scalh{H^{*}(\tau)R^{-1}(\tau)H(\tau)\zeta(\tau)}{\zeta(\tau)}d\tau\label{eq:3-19}
	\end{equation}
	If we call $x_{zt}(\tau)$ the optimal state and define $q_{zt}(\tau)$ as in (\ref{eq:3-16}) we obtain by standard methods that the optimal
	controls $\hat{\xi}$ and $\hat{v}(\cdot)$ are given by 
	\begin{equation}
		\hat{\xi}=-\Pi_{0}^{\frac{1}{2}}q_{zt}(t_0),\;\hat{v}(\tau)=-Q^{\frac{1}{2}}(\tau)G^{*}(\tau)q_{zt}(\tau)\label{eq:3-20}
	\end{equation}
	and the pair $x_{zt}(\tau),q_{zt}(\tau)$ solves (\ref{eq:3-16}). 
	
	\section{IDENTIFYING THE RKHS}\label{sec:rkhs}
	
	We can draw a parallel between the above reasoning of \Cref{sec:impulse-response} and going from a partial differential equation (PDE) to its calculus of variations formulation. Similarly, to solve a linear PDE one can look for its Green kernel, a.k.a.\ its fundamental solution. This is precisely what we have done in \eqref{sec:twoPtsBdry} for \eqref{eq:twoPtsBdry}. Now we will take advantage of \eqref{eq:3-19} to identify the spaces of functions we are optimizing over. This is the same thing as identifying the reproducing kernel Hilbert spaces of the kernels of Theorem~\ref{thm:main_noyaux}. In general it is a hard task given a PDE to find its Green kernel and the vector space of functions where the solutions live \citep[Chapter 1]{saitoh16theory}. But this task will be made much simpler owing to the fact that we are not looking at any Hilbert space, but at an RKHS:
	\begin{definition}\label{def_vRKHS}
		Let $\T$ be a non-empty set. A Hilbert space $(\Hk(\T),\left<\cdot,\cdot\right>_{K})$ of $\R^n$-vector-valued functions defined on $\T$ is called a vRKHS if there exists a matrix-valued kernel $K_\T:\T \times \T \rightarrow \cL(\R^{n,*},\R^n)$ such that the \emph{reproducing property} holds: for all $t \in \T,\, \b p\in\R^{n,*} $, we have $K_\T(\cdot,t)\b p \in \Hk(\T)$ and for all $\b f \in \Hk(\T)$, $\scalh{\b p}{\b f(t)}= \left<\b f,K_\T(\cdot,t)\b p\right>_{K}$.
	\end{definition}
	\noindent \tb{Remark:} It is well-known that by Riesz's theorem, an equivalent definition of a vRKHS is that, for every $t \in \T$ and $\b p\in\R^n$, the evaluation functional $\b f\in\Hk(\T) \mapsto \scalh{\b p}{\b f(t)}\in\RR$ is continuous. There is also a one-to-one correspondence between the kernel $K_\T$ and the vRKHS $(\Hk(\T),\left<\cdot,\cdot\right>_{K})$ \citep[see e.g.][Proposition 2.3]{carmeli06vector}. Moreover, by symmetry of the scalar product, the matrix-valued kernel has a Hermitian symmetry, i.e.\ $K_\T(s,t)=K_\T(t,s)^*$ for any $s,t\in\T$. We refer to \citet{carmeli10vector} and references therein for more on this topic of operator-valued kernels.
	
	
	\subsection{$K$: A PRIMAL RKHS OF TRAJECTORIES}
	
	We want to show that with the kernel $K(s,t | T)$ defined in (\ref{eq:3-3}), 
	\begin{align}
		K(s,t | T)&=\Phi_{F,\Sigma}(s,t_0)\Pi_{0}^{\frac{1}{2}}(\Id+\Pi_{0}^{\frac{1}{2}}\Sigma(t_0)\Pi_{0}^{\frac{1}{2}})^{-1}\Pi_{0}^{\frac{1}{2}}\Phi_{F,\Sigma}^{*}(t,t_0)+\int_{t_0}^{\min(s,t)}\Phi_{F,\Sigma}(s,\tau)G(\tau)Q(\tau)G^{*}(\tau)\Phi_{F,\Sigma}^{*}(t,\tau)d\tau\label{eq:kerneK}
	\end{align}
	we can associate a Hilbert space, for which it is the reproducing kernel. As intuited from \eqref{eq:3-18}, it is a space of functions in $H^{1}(t_0,T;\R^{n})$ defined as follows
	\begin{multline}
		\Sx=\{x(\cdot)\in H^{1}(t_0,T;\R^{n})|\exists\:u(\cdot)\in L^{2}(t_0,T;\R^{p}),\xi\in \R^{n}\:\text{ s.t.\ }\\ \dfrac{d}{d\tau}x=F(\tau)x(\tau)+G(\tau)Q^{\frac{1}{2}}(\tau)u(\tau),x(t_0)=\Pi_{0}^{\frac{1}{2}}\xi\}.\label{eq:4-1}
	\end{multline}
	Recall that we only required $R(\cdot)$ to be invertible, improving over \citet{aubin2022operator}. For a given trajectory $x(\cdot)$, in case there are several pairs $(u(\cdot),\xi)$ which satisfy the relations (\ref{eq:4-1}) we call the representative
	of $x(\cdot)$ the pair with minimal norm $\int_{t_0}^{T}\|u(t)\|^{2}dt +\|\xi\|^{2}$. It is uniquely defined. We can identify this representative as the unique pair $(u(\cdot),\xi)$ satisfying (\ref{eq:4-1}) for a given $x(\cdot)$ and 
	\begin{equation}
		\int_{t_0}^{T} \scalh{u(s)}{\widetilde{u}(s)}ds=0,\scalh{\xi}{\widetilde{\xi}}=0,\forall\, \widetilde{u}(\cdot),\widetilde{\xi} \text{ such that \:}G(s)Q^{\frac{1}{2}}(s)\widetilde{u}(\cdot)=0,\text{for a.e.}\,s\in(t_0,T),\:\Pi_{0}^{\frac{1}{2}}\widetilde{\xi}=0\label{eq:4-100}
	\end{equation}
	The vector space $\Sx$ is a subspace of the Sobolev space $H^{1}(t_0,T;\R^{n})$. We equip
	$\Sx$ with a scalar product derived from \eqref{eq:3-19}
	\begin{equation}
		\scalidx{x^{1}(\cdot)}{x^{2}(\cdot)}{\Sx}=\scalh{\xi^{1}}{\xi^{2}}+\scalh{\Sigma_Tx^{1}(T)}{x^{2}(T)}+\int_{t_0}^{T}\scalh{u^{1}(s)}{u^{2}(s)}ds+\int_{t_0}^{T}\scalh{H^{*}(s)R^{-1}(s)H(s)x^{1}(s)}{x^{2}(s)}ds\label{eq:4-2}
	\end{equation}
	where $(\xi^{1},u^{1}(\cdot))$; $(\xi^{2},u^{2}(\cdot))$ are the representatives
	of $x^{1}(\cdot),x^{2}(\cdot)$ respectively. 
	
	
	\begin{thm}
		\label{theo4-1}The vector space $(\Sx,\scalidx{\cdot}{\cdot}{\Sx})$ is a reproducing
		kernel Hilbert space, and $K(s,t|T)$ defined by (\ref{eq:3-3}) is
		its corresponding reproducing kernel.
	\end{thm}
	\begin{proof}
		The vector space $(\Sx,\scalidx{\cdot}{\cdot}{\Sx})$ is clearly a pre-complete Hilbert space by bilinearity of its inner product. It is furthermore complete since every Cauchy sequence converges within. Indeed, consider a Cauchy sequence $(x^{k}(\cdot))_{k\in\N}$ in $\Sx$, namely one for which there exists $(u^{k}(\cdot),\xi^{k})$, which is the representative
		of $x^{k}(\cdot)$, satisfying
		\[
		\dfrac{d}{ds}x^{k}=F(s)x^{k}(s)+G(s)Q^{\frac{1}{2}}(s)u^{k}(s), \quad x^{k}(t_0)=\Pi_{0}^{\frac{1}{2}}\xi^{k}
		\]
		and
		\begin{multline*}
			\|\xi^{k}-\xi^{l}\|^{2}+\|\Sigma_T^{\frac{1}{2}}(x^{k}(T)-x^{l}(T))\|^{2}+\int_{t_0}^{T}\|u^{k}(s)-u^{l}(s)\|^{2}ds\\+\int_{t_0}^{T}\scalh{H^{*}(s)R^{-1}(s)H(s)(x^{k}(s)-x^{l}(s))}{(x^{k}(s)-x^{l}(s)}ds\xrightarrow[k,l\rightarrow+\infty]{} 0.
		\end{multline*}
		Consequently $\xi^{k}$$, u^{k}(\cdot)$ are Cauchy sequences in $\R^{n},L^{2}(t_0,T;\R^{m})$ respectively, so $\xi^{k}\rightarrow\xi,u^{k}(\cdot)\rightarrow u(\cdot)$
		and necessarily $x^{k}(\cdot)\rightarrow x(\cdot)$ in $H^{1}(t_0,T;\R^{n})$ with 
		\[
		\dfrac{d}{ds}x=F(s)x(s)+G(s)Q^{\frac{1}{2}}(s)u(s),x(t_0)=\Pi_{0}^{\frac{1}{2}}\xi.
		\]
		Since $(\xi^{k},u^{k}(\cdot))$ are the representatives of $x^{k}(\cdot)$ they
		satisfy  $\int_{t_0}^{T}  \scalh{u^{k}(s)}{\widetilde{u}(s)}ds=0$, $ \scalh{\xi^{k}}{\widetilde{\xi}}=0$
		for any pair $(\widetilde{u}(\cdot),\widetilde{\xi})$ $\text{such that \:G(s)\ensuremath{Q^{\frac{1}{2}}}(s)}\widetilde{u}(s)=0,\text{a.e.}\,s\in(t_0,T),\:\Pi_{0}^{\frac{1}{2}}\widetilde{\xi}=0$.
		But then $(u(\cdot),\xi)$ satisfies also $\int_{t_0}^{T}\scalh{u(s)}{\widetilde{u}(s)}ds=0, \scalh{\xi}{\widetilde{\xi}}=0$, 
		which implies that $(u(\cdot),\xi)$ is the representative of $x(\cdot)$. Hence $x^{n}(\cdot)\rightarrow x(\cdot)$ in the sense of the norm of $\Sx$ and $\Sx$ is a Hilbert space. We now want to show that $\Sx$ is a RKHS and that its corresponding kernel is $K(s,t | T)$. We must check two facts as per Definition~\ref{def_vRKHS}
		\begin{gather}
			K(\cdot,t | T)z\in\Sx\label{eq:4-4}\\
			\scalidx{x(\cdot)}{K(\cdot,t | T)z}{\Sx}=\scalh{x(t)}{z},\forall x(\cdot)\in\Sx,t\in[t_0,T],z\in \R^{n,*}\label{eq:4-5}
		\end{gather}
		By definition, $K(\cdot,t|T)z=x_{zt}(\cdot)$. Now, from formulas (\ref{eq:3-16})
		$x_{zt}(\cdot)$ belongs to $\Sx$,  with 
		\begin{equation}
			v_{zt}(\tau)=-Q^{\frac{1}{2}}(\tau)G^{*}(\tau)q_{zt}(\tau)+Q^{\frac{1}{2}}(\tau)G^{*}(\tau)\chi_{zt}(\tau)\1_{\tau<t},\quad \xi_{zt}=-\Pi_{0}^{\frac{1}{2}}q_{zt}(t_0)+\Pi_{0}^{\frac{1}{2}}\chi_{zt}(t_0)\label{eq:4-6}
		\end{equation}
		We claim that $(v_{zt}(\cdot),\xi_{zt})$ is the representative of $x_{zt}(\cdot)$.
		Indeed the pair $(v_{zt}(\cdot), \xi_{zt})$ satisfies \eqref{eq:4-100} by using a transposition. 
		It remains to check (\ref{eq:4-5}). Consider $x(\cdot)\in\Sx$, such that 
		\[
		\dfrac{d}{d\tau}x=F(\tau)x(\tau)+G(\tau)Q^{\frac{1}{2}}(\tau)u(\tau),x(t_0)=\Pi_{0}^{\frac{1}{2}}\xi
		\]
		assuming that $(u(\cdot),\xi)$ is the representative. We have 
		\begin{align*}
			&\scalidx{x(\cdot)}{x_{zt}(\cdot)}{\Sx}=\scalh{\xi}{(-\Pi_{0}^{\frac{1}{2}}q_{zt}(t_0)+\Pi_{0}^{\frac{1}{2}}\chi_{zt}(t_0))}+\scalh{\Sigma_Tx(T)}{x_{zt}(T)}\\
			&+\int_{t_0}^{T} \scalh{u(\tau)}{-Q^{\frac{1}{2}}(\tau)G^{*}(\tau)q_{zt}(\tau)+Q^{\frac{1}{2}}(\tau)G^{*}(\tau)\chi_{zt}(\tau)\1_{\tau<t} }d\tau+\int_{t_0}^{T} \scalh{H^{*}(\tau)R^{-1}(\tau)H(\tau)x(\tau)}{x_{zt}(\tau)}d\tau
		\end{align*}
		Using relations (\ref{eq:3-16}) and the equation for $x(\cdot)$ with
		integration by parts, it is easy to conclude that
		\begin{align*}
			\scalidx{x(\cdot)}{x_{zt}(\cdot)}{\Sx}&=\scalh{x(t_0)}{\chi_{zt}(t_0)-q_{zt}(t_0)}+\scalh{x(T)}{q_{zt}(T)}+\int_{t_0}^{T}\scalh{\dfrac{d}{d\tau}x(\tau)-F(\tau)x(\tau)}{\chi_{zt}(\tau)\1_{\tau<t}-q_{zt}(\tau)}d\tau\\
			&\hspace{2cm}+\int_{t_0}^{T}\scalh{x(\tau)}{-\dfrac{d}{d\tau}q_{zt}(\tau)-F^*(\tau)q_{zt}(\tau)}d\tau\\
			&\hspace{-2cm}=\left[\scalh{x(\tau)}{q_{zt}(\tau)}\right]^{t_0}_{\tau=T}+\scalh{x(T)}{q_{zt}(T)}+\int_{t_0}^{t}\scalh{\dfrac{d}{d\tau}x(\tau)}{\chi_{zt}(\tau)}d\tau+\int_{t_0}^{t}\scalh{x(\tau)}{\dfrac{d}{d\tau}\chi_{zt}(\tau)}]d\tau=\scalh{x(t)}{z}
		\end{align*}
		which proves $(\ref{eq:4-5})$ and concludes the proof of the theorem.
		
	\end{proof}
	
	\begin{remark}[No measurements and Gramian of controllability]
		In the simpler case where $\Sigma_{T}=0$, $H(\cdot)\equiv 0$, we have that $\Sigma(\cdot)\equiv 0$. For $\Pi_{0}=0$ and $Q(\cdot)\equiv \Id$, \eqref{eq:3-19} becomes $J_x(\xi,u(\cdot))=\int_{t_0}^{T}\|u(\tau)\|^{2}d\tau$ and we obtain
		\begin{equation*}
			K(T,T | T)=\int_{t_0}^{T}\Phi_{F}(T,\tau)G(\tau)Q(\tau)G^{*}(\tau)\Phi_{F}^{*}(T,\tau)d\tau
		\end{equation*}
		which is precisely the Gramian of controlability. Its straightforward relation with kernels for optimal control was further discussed in \citet[p8]{aubin2020hard_control}.
	\end{remark}
	\subsection{$\Lambda$: A DUAL RKHS OF INFORMATION VECTORS}\label{sec:dualRKHS}
	We have seen with formula (\ref{eq:2-17}) that the optimal operator
	$\hat{S}_s(t | T)$ is directly related to an optimal control problem (\ref{eq:2-13}),(\ref{eq:2-14}).
	This problem concerns a backward evolution. Similarly, formula (\ref{eq:2-200})
	has been related to a two-point boundary value problem (\ref{eq:twoPtsBdry}). This one is easily interpreted as the necessary and sufficient conditions of optimality of another control problem, namely 
	\begin{equation}
		-\dfrac{d}{dt}\lambda(t)=F^{*}(t)\lambda(t)+H^{*}(t)v(t), \quad \lambda(T)=\Sigma_T^{\frac{1}{2}} z\label{eq:5-8}
	\end{equation}
	and the payoff to minimize is a slight variation of (\ref{eq:2-14}), 
	\begin{equation}
		J_\lambda(z,v(\cdot))=\scalh{\Pi_{0}\lambda(t_0)}{\lambda(t_0)}+\norh{z}^2+\int_{t_0}^{T}\scalh{G(t)Q(t)G^{*}(t)\lambda(t)}{\lambda(t)}dt+\int_{t_0}^{T}\scalh{R(t)v(t)}{v(t)}dt.\label{eq:5-9}
	\end{equation}
	Recall that we did not require $\Sigma_T$ to be invertible but that $R(t)$ is invertible. These problems lead to the kernel
	\begin{align}
		\Lambda(s,t|T)&=\Phi_{F,\Pi}^{*}(T,s)\Sigma_T^{\frac{1}{2}}(\Id+\Sigma_T^{\frac{1}{2}}\Pi(T)\Sigma_T^{\frac{1}{2}})^{-1}\Sigma_T^{\frac{1}{2}}\Phi_{F,\Pi}(T,t)+\int_{\max(s,t)}^{T}\Phi_{F,\Pi}^{*}(\tau,s)H^{*}(\tau)R^{-1}(\tau)H(\tau)\Phi_{F,\Pi}(\tau,t)d\tau\label{eq:kerneLa}
	\end{align}
	and to a space of information vectors $\Sla$
	\begin{multline}
		\Sla=\{\lambda(\cdot)\in H^{1}(t_0,T;\R^{n})|\:v(\cdot)\in L^{2}(t_0,T;\R^{m}),z\in\R^n\:\text{ s.t.\ }\\ -\dfrac{d}{dt}\lambda(t)=F^{*}(t)\lambda(t)+H^{*}(t)v(t),\;\lambda(T)=\Sigma_T^{\frac{1}{2}} z\}.\label{eq:5-10}
	\end{multline}
	Again if there are several $(z,v(\cdot))$ for the same trajectory $\lambda(\cdot)$, we define the representative of $\lambda(\cdot)$ as the unique $(z,v(\cdot))$ satisfying (\ref{eq:5-10}) which minimizes $\int_{t_0}^{T}\scalh{R(t)v(t)}{v(t)}dt+\|z\|^{2}$ and, similarly to \eqref{eq:4-100}, we have
	\begin{equation}\label{eq:4-101}	
		\int_{t_0}^{T} \scalh{R(s)v(s)}{\widetilde{v}(s)}ds=0,\scalh{z}{\widetilde{z}}=0,\forall\, \widetilde{v}(\cdot),\widetilde{z} \text{ such that \:}H^*(s)\widetilde{v}(\cdot)=0,\text{for a.e.}\,s\in(t_0,T),\:\Sigma_{T}^{\frac{1}{2}}\widetilde{z}=0.
	\end{equation}
 We then equip $\Sla$ with the inner product derived from the quadratic form \eqref{eq:5-9}.
	\begin{thm}
		\label{theo5-2} The Hilbert space $(\Sla,\scalidx{\cdot}{\cdot}{\Sla})$ is a reproducing
		kernel Hilbert space and the corresponding kernel is $\Lambda(s,t|T)$. 
	\end{thm}
	
	\begin{proof} The proof is exactly the same as that of Theorem~\ref{theo4-1}, noticing that the RKHS $\Sx$ corresponds to $\Sla$ when considering the permutations $t\leftrightarrow T-(t-t_0)$, $F(t)\leftrightarrow F^*(t)$, $\Pi_0\leftrightarrow \Sigma_T$, $G(t)Q^{\frac{1}{2}}(t)\leftrightarrow H^*(t)R^{-\frac{1}{2}}(t)$.
	\end{proof}
	
	\begin{remark}[No plant noise and Gramian of observability]
		In the simpler case where $\Pi_{0}=0$, $G(\cdot)\equiv 0$, we have that $\Pi(\cdot)\equiv 0$ and $\Sx$ is finite-dimensional. For $\Sigma_{T}=0$ and $R(\cdot)\equiv \Id$, \eqref{eq:5-9} becomes $J_\lambda(z,v(\cdot))=\int_{t_0}^{T}\norh{v(t)}^2dt$ and we obtain
		\begin{equation*}
			\Lambda(t_0,t_0 | T)=\int_{t_0}^{T}\Phi_{F}^{*}(\tau,t_0)H^{*}(\tau)H(\tau)\Phi_{F}(\tau,t_0)d\tau
		\end{equation*}
		which is precisely the Gramian of observability. This setting is discussed in \citet[Example 16.3.3]{Kailath2000}. The covector $\lambda(\cdot)$, through its interpretation as a Lagrange multiplier as in \Cref{sec:detour_estimator}, is related to a notion of sensitivity of the solution $\hat{x}(\cdot)$ to the ``constraints'' $(y_t)_{t\in[t_0,T]}$, and thus here to the ability to recover $x_0$ from the observations.
	\end{remark}
	By Loève's theorem \citep[Theorem 27, p57]{berlinet04reproducing}, we know that $\Lambda(s,t)$ is the covariance of a stochastic process. It is clearly not forward Markovian. It is possible to introduce a backward Markovian Gaussian process which would allow to recover \eqref{eq:mvle-dual2}, with technicalities covered in \citet[Appendix Chap.\ 16]{Kailath2000}. We leave this identification to future developments.
	
	\subsection{$K$ AS ERROR COVARIANCE, FORMULAS RELATING $K$ AND $\Lambda$}
	We can now prove the claimed result \eqref{eq:kernel_def} stating that the kernel $K$ corresponds to the covariance of the error. We can also relate through a simple formula the two kernels $K$ and $\Lambda$.
	\begin{prop}\label{prop:covar_BF_formula} Let $\Sigma_{T}=0$ and recall that $\epsilon_{\hat S_s}(s | T)=x(s)-\hat{x}(s|T)$, then
		\begin{equation}
			K(s,t | T)=\bE[\epsilon_{\hat S_s}(s | T)(\epsilon_{\hat S_t}(t | T))^{*}]\label{eq:kernel_covar}
		\end{equation}
		Moreover, we have, for all $s,t\in[t_0,T]$,
		\begin{equation}
			K(s,t | T)=\Pi(s)\Phi_{F,\Pi}^{*}(t,s)\1_{s\le t}+\Phi_{F,\Pi}(s,t)\Pi(t)\1_{s> t}-\Pi(s)\Lambda(s,t|T)\Pi(t).\label{eq:BF_formula}
		\end{equation}
		In particular $K(T,T|T)=\Pi(T)$ and $	K(t,t | T)=\Pi(t)-\Pi(t)\Lambda(t,t|T)\Pi(t)$.
	\end{prop}
	Equation \eqref{eq:BF_formula} can be seen as an extension to $s\neq t$ of Bryson-Frazier formulas \citep[Theorem 16.5.1]{Kailath2000}. For $s=t$, it gives that $\bE[\epsilon_{\hat S_s}(s | T)(\epsilon_{\hat S_s}(s | T))^{*}] \preccurlyeq \Pi(s)= \bE[\epsilon_{\hat S_s}(s | s)(\epsilon_{\hat S_s}(s | s))^{*}] $. In other words, observations on the interval $[s,T]$ allow to lower the variance of the error at time $s$. For $t=T$ and $\Sigma_{T}=0$, $\Lambda(s,T|T)=0$ and we recover the result of \citet[Lemma 16.5.1]{Kailath2000}.
	
	\begin{proof}
		Since for  $l_\mu(\cdot)\equiv 0$ and $l_{\nu,s}(\tau)=-\bla \delta_s(\tau)$, the system \eqref{eq:2-15} with a jump condition corresponds to the two-point boundary system \eqref{eq:twoPtsBdry}. With the same computations that led to \eqref{eq:2-12}, using \eqref{eq:2-17} then integration by parts, we obtain that
		\begin{align*}
			&\scalh{\bE[\epsilon_{\hat S_s}(s | T)(\epsilon_{\hat S_s}(t | T))^{*}]\bla}{\bla}=\scalh{\Pi_{0}\lambda_s(t_0)}{\lambda_t(t_0)}+\int_{t_0}^{T}\scalh{G(\tau)Q(\tau)G^{*}(\tau)\lambda_s(\tau)}{\lambda_t(\tau)}d\tau\\
			&\hspace{6cm}+\int_{t_0}^{T}\scalh{H^{*}(\tau)R^{-1}(\tau)H(\tau)\hat{\gamma}_s(\tau)\bla}{\hat{\gamma}_t(\tau)}d\tau\\
			&\stackrel{\eqref{eq:2-15}}{=}-\scalh{\lambda_s(t_0)}{\hat{\gamma}_t(t_0)}+\int_{t_0}^{T}\scalh{\lambda_s(\tau)}{F(\tau)\hat{\gamma}_s(\tau)-\dfrac{d\hat{\gamma}_s}{du}}d\tau-\int_{t_0}^{T}\scalh{\dfrac{d\hat{\lambda}_s}{du}+F^{*}(\tau)\hat{\lambda}_s(\tau)}{\hat{\gamma}_t(\tau)}d\tau\\
			&=-\scalh{\lambda_s(t_0)}{\hat{\gamma}_t(t_0)}-[\scalh{\lambda_{s}(\tau)}{\hat{\gamma}_t(\tau)}]_{t_0}^s-[\scalh{\lambda_{s}(\tau)}{\hat{\gamma}_t(\tau)}]_{s^+}^T\stackrel{\eqref{eq:2-15}}{=}-\scalh{\bla}{\hat{\gamma}_t(s)}\stackrel{\eqref{eq:Wiener_K}}{=}\scalh{\bla}{K(s,t|T)\bla}.
		\end{align*}
		This proves \eqref{eq:kernel_covar}. For \eqref{eq:BF_formula}, define $r_t(\tau)=\hat{\gamma}_t(\tau)+\Pi(\tau)\lambda_{t}(\tau)$ as in \eqref{eq:error_edo}. Then similarly to \eqref{eq:r_explicite}-\eqref{eq:edo_nu_explicite}
		\begin{gather}
			r_t(\tau)=-\int_{t_{0}}^\tau \Phi_{F,\Pi}(\tau,\sigma)(\Pi(\sigma)\bla \delta_t(\sigma)) d\sigma=-\Phi_{F,\Pi}(\tau,t)\Pi(t)\1_{\tau> t}\bla. \label{eq:r_explicite2}\\
			-\dfrac{d}{d\tau}\lambda_{t}=(F^{*}(\tau)-H^{*}(\tau)R^{-1}(\tau)H(\tau))\lambda_{t}(\tau)+\bla \delta_t(\tau) -H^{*}(\tau)R^{-1}(\tau)H(\tau) \Phi_{F,\Pi}(\tau,t)\Pi(t)\1_{\tau> t}\bla.\label{eq:edo_nu_explicite2}
		\end{gather}
		Consequently, by the variation of constants formula,
		\begin{equation}
			\lambda_{t}(s)=\Phi_{F,\Pi}^{*}(t,s)\bla \1_{s\le t}-\int_{\max(s,t)}^{T}\Phi_{F,\Pi}^{*}(\tau,s)H^{*}(\tau)R^{-1}(\tau)H(\tau)\Phi_{F,\Pi}(\tau,t)\Pi(t)\bla d\tau.\label{eq:nu_explicite2}
		\end{equation}
		Since $\scalh{\bla}{K(s,t|T)\bla}=-\scalh{\bla}{\hat{\gamma}_t(s)}=\scalh{\bla}{\Pi(s)\lambda_{t}(s)-r_t(s)}$, \eqref{eq:r_explicite2}-\eqref{eq:nu_explicite2} yield \eqref{eq:BF_formula}.
		
	\end{proof}
	\begin{remark}[Relation with Fisher information]
		From \eqref{eq:kernel_covar}, we know that the kernel $K$ is the optimal error covariance for linear estimators. The Fisher information matrix of estimating $(x_t)_{t\in[t_0,T]}$ given $(y_t)_{t\in[t_0,T]}$ is well-defined if the unknown $(x_t)_{t\in[t_0,T]}$ is deterministic. In that case, since $(y_t)_{t\in[t_0,T]}$ is a Gaussian process and the measurements are linear, the integral operator with kernel $K$ thus saturates the Cramér-Rao inequality and coincides with the inverse of the Fisher information matrix. For random $(x_t)_{t\in[t_0,T]}$, defining properly the Bayesian posterior analogue of the Cramér-Rao lower bound \citep{VanTrees_1968} in this infinite-dimensional context would require an undue level of generality, as hinted at in \Cref{sec:detour_estimator}. However, for a deterministic $(x_t)_{t\in[t_0,T]}$, we have $G(\cdot)\equiv 0$ and $\Pi_0=0$ as considered previously. The only parameter of $(x_t)_{t\in[t_0,T]}$ to identify is $x_0$, so the Fisher information matrix writes as $\int_{t_0}^{T} \frac{\partial \bar y(t)}{\partial x_0}\frac{\partial \bar y(t)}{\partial x_0}^* dt$ which through simple calculations is shown to be equal to the Gramian of observability (already noticed in \citet[Remark (h)]{Kalman1961NewRI}). The fact that the Fisher information defines a kernel was unveiled by \citet{Jaakkola98exploitinggenerative}. Since the posterior Cramér-Rao lower bound of our estimation problem does not seem to have been written yet, much remains to be done in this direction of research to draw connections with information theory and with the score vector. Applications would already exist in sensor placement and D-optimal design.
	\end{remark}
	\subsection{DUAL DETERMINISTIC PROBLEMS ON RKHSs}\label{sec:dual_problem_rkhs}
	
	So far, we identified the RKHSs associated with the two Green functions $K$ and $\Lambda$. However the primary modern use of reproducing kernels is as convenient spaces for optimization problems such as the GP regression \eqref{eq:gpr} highlighted earlier on. This is exemplified in the following formalization of \eqref{eq:mvle-ls}. Since we do not require any invertibility except for $R(t)$, we will use pseudo-inverses w.r.t.\ the Euclidean norm for the other operators, such as $\Pi_0^\ominus$, but this does not change the formalism. For $\widetilde{x}(\cdot)\in\Sx$, after making the covariances explicit, \eqref{eq:mvle-ls} writes as the traditional least-squares problem:
	\begin{multline}\label{eq:lse-kernel}
		L_x(\widetilde{x}(\cdot)):=\int_{t_{0}}^T \noridx{\widetilde{y}(t)- H(t)\widetilde{x}(t)}{R(t)^{-1}}^2dt+\noridx{G(t)^{\ominus}\left(\frac{d}{dt}\widetilde{x}-F(t)\widetilde{x}(t)\right)}{Q(t)^{\ominus}}^2dt+\scalh{\Pi_0^\ominus \widetilde{x}(t_0)}{\widetilde{x}(t_0)}+\scalh{\Sigma_T\widetilde{x}(T)}{\widetilde{x}(T)}\\ = \int_{t_{0}}^T \noridx{\widetilde{y}(t)- H(t)\widetilde{x}(t)}{R(t)^{-1}}^2dt+ \noridx{\widetilde{x}(\cdot)}{\Sx}^2-\int_{t_{0}}^T \noridx{ H(t)\widetilde{x}(t)}{R(t)^{-1}}^2dt,
	\end{multline}
	where we refer to footnote~\ref{foo:observation_continue} for the interpretation of $\widetilde{y}(t)$. The space $\Sx$ is actually the set of trajectories such that $L_x$ is finite. The term in $\Sigma_T$ can be seen as imposing a Gaussian prior on the terminal point.
	\begin{prop}\label{prop:primal_prob_rkhs}
		For a realization $\widetilde{y}(\cdot)$, the zero-mean part of the solution to the smoothing problem is given by
		\begin{multline}\label{eq:primal_prob_rkhs}
			\int_{t_0}^{T}K(\cdot,t | T)H^{*}(t)R^{-1}(t)\widetilde{y}(t)dt \\= \argmin_{\tilde{x}(\cdot)\in\Sx}L_x(\widetilde{x}(\cdot))= \noridx{R(t)^{\nicefrac{-1}{
				2}}\widetilde y(\cdot)}{L^2}^2+\noridx{\widetilde{x}(\cdot)}{\Sx}^2-2\scalidx{H^{*}(\cdot)R^{-1}(\cdot)\widetilde y(\cdot)}{\widetilde{x}(\cdot)}{L^2([t_0,T])}
		\end{multline}
	\end{prop}
	\begin{proof}
	By Corollary~\ref{prop3-1}, we know that the l.h.s.\ is the solution of the smoothing problem. We just have to prove the equality. Set $\T=[t_0,T]$. By our assumptions on $F,G,Q,H$, we obtain that $K(\cdot,\cdot|T)\in L^\infty(\T\times\T,\cL(\R^{n,*},\R^n))$. Consequently the following kernel integral operator $\scK:L^2(\T,\R^{n,*})\rightarrow L^2(\T,\R^{n})$ is self-adjoint and bounded \citep[Proposition 4.2]{carmeli06vector},
	\begin{align*}
	 	(\scK f) (s):= \int_{t_0}^T K(s,t|T)f(t)dt \text{ for $f\in L^2(\T,\R^{n,*})$}.
	\end{align*}
	Furthermore, by the reproducing property \eqref{eq:4-5}, we have that
	\begin{align*}
	 \forall f\in L^2(\T,\R^{n,*}), \, g\in\Sx,\;	\scalidx{g}{f}{L^2(\T)}=\int_{t_{0}}^T \scalh{f(t)}{g(t)}dt=\scalidx{g}{\int_{t_0}^T K(\cdot,t|T)f(t)dt}{\Sx}=\scalidx{g}{\scK f}{\Sx}.
	\end{align*}
	 Consequently \eqref{eq:primal_prob_rkhs} boils down to minimizing the strongly convex function $\noridx{\widetilde{x}(\cdot)}{\Sx}^2-2\scalidx{\scK [H^{*}R^{-1}\widetilde y(\cdot)]}{\widetilde{x}}{\Sx}$, which gives $\hat{x}(\cdot)=\scK [H^{*}R^{-1}\widetilde y(\cdot)]$ as expected. This concludes the proof.
	\end{proof}
	\begin{remark}[Analogy with Sobolev spaces in calculus of variations]
		The above problem \eqref{eq:primal_prob_rkhs} may seem unfamiliar to kernel practitioners. Indeed it would be a kernel ridge regression if not for the last term in \eqref{eq:lse-kernel}, as it is not a quadratic data fitting term with a quadratic regularizer, but a linear data fitting term. Nevertheless it precisely matches the variational formulation ($\int \norh{\nabla u(t)}^2 dt+\scalidx{f}{u}{L^2}$) of a PDE like Poisson's equation ($\Delta u= f$) with null boundary condition. The emphasis we have put on the Sobolev-like space $\Sx$ over which we optimize is no different than the focus on the Sobolev space $H^1_0$ to study Poisson's equation. \citet{Rockafellar1987} treated Linear-Quadratic optimal control as a form of quadratic programming. The Hilbert spaces of trajectories and information vectors we defined are complementary to his discussion.\footnote{\citet{Rockafellar1987} takes the control in $L^1$ and considers additional constraints, making it closer to optimization problems over Banach spaces $W^{1,\infty}$ or $W^{1,1}$ which we do not cover in our Hilbertian setting.} Considering a dual convex problem to an optimal control one is admittedly less frequent in the control community, even though it may have some computational advantages \citep[see][and references therein]{burachik2014duality}. Here we derive the dual problem to \eqref{eq:primal_prob_rkhs}, which corresponds formally to \eqref{eq:mvle-3}.
	\end{remark}
	\begin{prop}\label{prop:dual_prob_rkhs}
		Decompose $R^{-1}(t)\widetilde y(t)$ into the sum of two vectors $\proj\nolimits^{\|\cdot\|_{R(t)}}_{\Img H(t)}(\widetilde{y}(t))\in R(t)^{-1}\Img H(t)$ and $\proj\nolimits^{\|\cdot\|_{R(t)}}_{\Ker H^*(t)}(\widetilde{y}(t))\in \Ker H^*(t)$. The convex dual problem to \eqref{eq:primal_prob_rkhs} is, for the adjoint control $v(\cdot)$ defined as in \eqref{eq:5-10}, 
		\begin{multline}\label{eq:dual_prob_rkhs}
			\min_{\lambda(\cdot)\in\Sla}L_\lambda(\lambda(\cdot))=\noridx{\lambda(\cdot)}{\Sla}^2 - 2\int_{t_{0}}^T \scalh{R(t)\proj\nolimits^{\|\cdot\|_{R(t)}}_{\Img H(t)}(\widetilde{y}(t))}{v(t)}dt-\noridx{R(\cdot)^{\nicefrac{1}{
						2}}\proj\nolimits^{\|\cdot\|_{R(\cdot)}}_{\Ker H^*(\cdot)}(\widetilde{y}(\cdot))}{L^2}^2\\
			=\int_{t_{0}}^T \noridx{\proj\nolimits^{\|\cdot\|_{R(t)}}_{\Img H(t)}(\widetilde{y}(t))- v(t)}{R(t)}^2dt+\int_{t_0}^{T}\scalh{G(t)Q(t)G^{*}(t)\lambda(t)}{\lambda(t)}dt  +\scalh{\Pi_{0}\lambda(t_0)}{\lambda(t_0)}\\
			+\scalh{\Sigma_T^\ominus \lambda(T)}{\lambda(T)}-\noridx{R(\cdot)^{\nicefrac{-1}{
							2}}\widetilde{y}(\cdot)}{L^2}^2.
		\end{multline}
	\end{prop}
	\begin{proof}
		For any given $\lambda(\cdot)\in C^1(\T,\R^{n,*})$ and $\lambda_0\in\R^{n,*}$, we introduce the Lagrangian $L$ over $x(\cdot)\in C^1(\T,\R^{n})$, $u(\cdot)\in L^2(\T,\R^{d})$, $\xi\in\R^n$,
		\begin{multline*}
			L_{tot}(\xi, u(\cdot),x(\cdot),\lambda(\cdot),\lambda_0)=\|\xi\|^{2}+\scalh{\Sigma_Tx(T)}{x(T)}+\int_{t_0}^{T}\|u(\tau)\|^{2}d\tau+\int_{t_0}^{T}\scalh{H^{*}(\tau)R^{-1}(\tau)(H(\tau)x(\tau)-2 \widetilde y(\tau))}{x(\tau)}d\tau\\ + 2\int_{t_{0}}^T \scalh{\lambda(t)}{\dfrac{d}{d\tau}x-F(\tau)x(\tau)-G(\tau)Q^{\frac{1}{2}}(\tau)u(\tau)}d\tau+2\scalh{x(t_0)-\Pi_{0}^{\frac{1}{2}}\xi}{\lambda_0}+\noridx{R(\cdot)^{\nicefrac{-1}{
						2}}\widetilde y(\cdot)}{L^2}^2.
		\end{multline*}
		Integrating by parts, and minimizing over $(\xi, x(\cdot), u(\cdot), x(T), x(t_0))$, and we obtain that
		\begin{align*}
			&\partial_\xi &\xi=\Pi_{0}^{\frac{1}{2}} \lambda_0, \quad \quad	&\partial_{u(\cdot)}  &u=Q^{\frac{1}{2}} G^* \lambda,\\
			&\partial_{x(T)} &\lambda(T)=-\Sigma_T x(T), \quad \quad
			&\partial_{x(t_0)} &\lambda_0=\lambda(t_0),\\
			&\partial_{x(\cdot)} &0=\dfrac{d}{dt}\lambda(t)+F^{*}(t)\lambda(t)-H^{*}(t)R^{-1}(t)(H(t) x(t)-\widetilde y(t)),
		\end{align*}
	whence $\lambda(\cdot)\in\Sla$. Decompose $R^{-1}(t)\widetilde y(t)$ into $v_{H,y}(t)+v_{H*,y}(t)$ with $R(t)v_{H,y}(t) \in \Img H(t)$ and $v_{H*,y}(t) \in \Ker H^*(t)$.  Denoting by $\hat x$ the optimum of $x(\cdot)$ and setting $v(t)=-R^{-1}(t)H(t)\hat x(t)+v_{H,y}(t)$, $z=-\Sigma_T \hat x(T)$, we obtain that  $v(t)$ is the representative of $\lambda(\cdot)\in\Sla$, since $H^{*}(t)R^{-1}(t)\widetilde y(t)=H^{*}(t)v_{H,y}(t)$ and for all $\widetilde{v}$ satisfying \eqref{eq:5-10} and $H^*(t)\widetilde{v}(t)=0$, we have obviously $\scalh{R(t)v(t)}{\widetilde{v}(t)}=0$, so \eqref{eq:4-101} holds, and for $z=-\Sigma_T \hat x(T)$ as well, showing $(v_{H,y}(\cdot,z))$ is the representative of $\lambda(\cdot)$. Consequently
	\begin{align*}
		L(\hat x(\cdot),\lambda(\cdot))&=-\scalh{\Pi_{0}\lambda(t_0)}{\lambda(t_0)}-\norh{z}^2-\int_{t_0}^{T}\scalh{G(t)Q(t)G^{*}(t)\lambda(t)}{\lambda(t)}dt\\
		&-\int_{t_0}^{T}\scalh{H^{*}(\tau)R^{-1}(\tau)(H(\tau)\hat x(\tau)}{\hat x(\tau)}d\tau+\noridx{R(\cdot)^{\nicefrac{-1}{
					2}}\widetilde y(\cdot)}{L^2}^2,\\
		\intertext{since $\int_{t_0}^{T}\scalh{H^{*}(\tau)R^{-1}(\tau)(H(\tau)\hat x(\tau)}{\hat x(\tau)}d\tau=\scalidx{R(v-v_{H,y})}{v-v_{H,y}}{L^2}$,}
			L(\hat x(\cdot),\lambda(\cdot))&=-\noridx{\lambda(\cdot)}{\Sla}^2 + 2\int_{t_{0}}^T \scalh{R(t)v_{H,y}(t)}{v(t)}dt-\noridx{R(\cdot)^{\nicefrac{1}{
						2}}\widetilde v_{H,y}(\cdot)}{L^2}^2+\noridx{R(\cdot)^{\nicefrac{-1}{
						2}}\widetilde y(\cdot)}{L^2}^2.
	\end{align*}
 	Since, by orthogonality $\noridx{R(\cdot)^{\nicefrac{-1}{
 				2}}\widetilde y(\cdot)}{L^2}^2=\noridx{R(\cdot)^{\nicefrac{1}{
 				2}}\widetilde v_{H,y}(\cdot)}{L^2}^2+\noridx{R(\cdot)^{\nicefrac{1}{
 				2}}\widetilde v_{H*,y}(\cdot)}{L^2}^2$, this concludes the proof.

	\end{proof}
	The two problems we obtained \eqref{eq:primal_prob_rkhs}-\eqref{eq:dual_prob_rkhs} are the proper formalization of \eqref{eq:mvle-ls}-\eqref{eq:mvle-3} when the function spaces are properly defined. Not using the RKHSs would be similar to not defining the Sobolev spaces when doing variational analysis of quadratic energies. Furthermore, unlike \eqref{eq:mvle-ls} and \eqref{eq:mvle-3}, which stem from stochastic or Bayesian viewpoints as presented in \Cref{sec:detour_estimator}, \eqref{eq:primal_prob_rkhs} and \eqref{eq:dual_prob_rkhs} can be easily modified to account for further side information, such as the sign of the trajectory to be reconstructed, which is a form of state constraint, or other objective functions. We can then use the full machinery of kernel representer theorems to derive a numerical solution, as in \citet{aubin2020hard_control}.

	\section{COMPLEMENTS ON THE FILTER AND SMOOTHER}\label{sec:filter}
	In this section, for completeness, we derive the expression of the Kalman filter in continuous time and the expression of the smoother based on  the innovation process, both being well-known. We take $\Sigma_{T}=0$ and $l_\nu(t)=H^{*}(t)R^{-1}(t)g(t)$ with $g(\cdot)\in L^{2}(t_0,T;\R^{m})$, and later we posit $g(t)dt=d\widetilde{y}(t)$. Going back to the expression (\ref{eq:2-8}) of the best estimate, and to the expression (\ref{eq:3-4}) of the operator $\widehat{S}_s(t | T))$, we can write the formula
	\begin{equation}
		\hat{x}(s | T)=\bar{x}(s)+\int_{t_0}^{T}K(s,t | T)H^{*}(t)R^{-1}(t)d\widetilde{y}(t),\label{eq:5-1}
	\end{equation}
	which gives for $s=T$ the Kalman filter formula 
	\begin{equation}
		\hat{x}(T|T)=\bar{x}(T)+\int_{t_0}^{T}K(T|T,t)H^{*}(t)R^{-1}(t)d\widetilde{y}(t).\label{eq:5-100}
	\end{equation}
	We want to show that the smoother at $s$ can be recursively expressed in terms of the sequence of filters between $s$ and $T$. As we deal with Linear-Quadratic-Gaussian estimation problems, we naturally find ourselves in a very thoroughly explored field and some expressions below are well-known, and we point out to where they could be found in textbooks \citep{Kailath2000,Bensoussan2018}. General kernel formulas are in any case new, and so is the connection of Kalman filtering and smoothing with the explicit RKHSs discussed in \Cref{sec:rkhs}.
	
	\subsection{REDERIVING THE KALMAN FILTER}

	Recall that we defined $r(t)=\hat{\mu}(t)+\Pi(t)\hat{\nu}(t)$ in \eqref{eq:error_edo}. By \eqref{eq:r_explicite},
	\begin{equation}
		\dfrac{dr}{dt}=(F(t)-\Pi(t)H^{*}(t)R^{-1}(t)H(t))r(t)+\Pi(t)H^{*}(t)R^{-1}(t)g(t),\quad r(t_0)=0.\label{eq:5-114}
	\end{equation}
	The importance of $r(t)$ is that, like $\Pi(t)$, it satisfies a forward differential equation and thus it does not depend on $T$.
	On the other hand, since $\hat{\nu}(T)=0$, $\widehat{\mu}(T)=r(T)$. But then from \eqref{eq:r_explicite} we can write 
	\begin{equation}
		r(T)=\int_{t_0}^{T}K(T|T,t)H^{*}(t)R^{-1}(t)g(t)dt=\int_{t_0}^{T}\widehat{S}(T|T,t)g(t)dt.\label{eq:5-115}
	\end{equation}
	If in (\ref{eq:5-115}) we substitute $g(t)dt$ by $d\widetilde{y}(t)$
	we obtain the stochastic differential equation (SDE) as given in \citet[Chapter 7]{Bensoussan2018}
	\begin{equation}
		dr=(F(t)-\Pi(t)H^{*}(t)R^{-1}(t)H(t))r(t)dt+\Pi(t)H^{*}(t)R^{-1}(t)d\widetilde{y}(t),\quad r(t_0)=0.\label{eq:5-116}
	\end{equation}
	and 
	\begin{equation}
		r(T)=\int_{t_0}^{T}\widehat{S}(T|T,t)d\widetilde{y}(t).\label{eq:5-120}
	\end{equation}
	From (\ref{eq:5-100}) it follows that the Kalman filter over $[0,t]$, i.e.\ only filtering without smoothing, denoted to simplify notation $\widehat{x}(t)=\widehat{x}(t|t)=\overline{x}(t)+r(t)$, thus satisfies the SDE
	\[
	d\widehat{x}(t)=(F(t)\widehat{x}(t)+f(t))dt+\Pi(t)H^{*}(t)R^{-1}(t)(d\widetilde{y}(t)-H(t)r(t)dt).
	\]
	Finally, we obtain the classical Kalman filter equation
	\begin{equation}
		d\widehat{x}(t)=(F(t)\widehat{x}(t)+f(t))dt+\Pi(t)H^{*}(t)R^{-1}(t)(dy(t)-(H(t)\widehat{x}(t)+h(t))), \quad \widehat{x}(t_0)=x_{0}.\label{eq:5-121}
	\end{equation}
	
	\subsection{EXPRESSION OF THE SMOOTHER IN TERMS OF INNOVATION}
	From \eqref{eq:edo_nu_explicite}, we know that
	\begin{equation}
		-\dfrac{d}{d\tau}\widehat{\nu}(\tau)=(F^{*}(\tau)-H^{*}(\tau)R^{-1}(\tau)H(\tau)\Pi(\tau))\widehat{\nu}(\tau)-H^{*}(\tau)R^{-1}(\tau)(g(\tau)-H(\tau)r(\tau)), \quad \widehat{\nu}(T)=0\label{eq:5-125}
	\end{equation}
	and thus, by the variation of constants formula,
	\begin{equation}
		\widehat{\nu}(s)=-\int_{s}^{T}\Phi_{F,\Pi}^{*}(t,s)H^{*}(t)R^{-1}(t)(g(t)-H(t)r(t))dt.\label{eq:5-126}
	\end{equation}
	Therefore, we have 
	\begin{equation}
		\widehat{\mu}(s)=r(s)+\Pi(s)\int_{s}^{T}\Phi_{F,\Pi}^{*}(t,s)H^{*}(t)R^{-1}(t)(g(t)-H(t)r(t))dt.\label{eq:5-127}
	\end{equation}
	If we replace in the equation $g(t)dt$ by $d\widetilde{y}(t)$ we
	must interpret $r(t)$ as the solution of the SDE (\ref{eq:5-116})
	. We then have 
	\begin{equation}
		\widehat{x}(s | T)=\overline{x}(s)+\widehat{\mu}(s)\label{eq:5-128}
	\end{equation}
	\[
	\widehat{x}(s)=\overline{x}(s)+r(s).
	\]
	From \eqref{eq:5-127}, we deduce that
	\begin{equation}
		\widehat{x}(s | T)=\widehat{x}(s)+\Pi(s)\int_{s}^{T}\Phi_{F,\Pi}^{*}(t,s)H^{*}(t)R^{-1}(t)(dy(t)-(H(t)\widehat{x}(t)+h(t))dt).\label{eq:5-129}
	\end{equation}
	Define the innovation process $e(\cdot)$ as follows
	\begin{equation}
		e(t)=y(t)-\int_{t_0}^{t}(H(\tau)\widehat{x}(\tau)+h(\tau))d\tau.\label{eq:5-130}
	\end{equation}
	It is well-known that the innovation process is a $\mathcal{Y}^{t}$ Wiener process with covariance matrix $R(t)$ \citep[see e.g.][Lemma 7.1]{Bensoussan2018}. In a nutshell we obtained that
	\begin{prop}
		\label{prop5-1} The Kalman smoother can be written as follows
		\begin{equation}
			\widehat{x}(s | T)=\widehat{x}(s)+\Pi(s)\int_{s}^{T}\Phi_{F,\Pi}^{*}(t,s)H^{*}(t)R^{-1}(t)de(t)\label{eq:5-200}
		\end{equation}
	\end{prop}
	This result is well-known and can be found in \citet[Lemma 16.5.1]{Kailath2000}. Consequently, the error of the Kalman smoother is given by 
	\begin{align}
		\widehat{\epsilon}(s | T)&=x(s)-\widehat{x}(s)-\Pi(s)\int_{s}^{T}\Phi_{F,\Pi}^{*}(t,s)H^{*}(t)R^{-1}(t)de(t)=\widehat{\epsilon}(s|s)-\Pi(s)\int_{s}^{T}\Phi_{F,\Pi}^{*}(t,s)H^{*}(t)R^{-1}(t)de(t)\label{eq:5-201}
	\end{align}
	It is obvious that, for $s<t$, the random variables $\widehat{\epsilon}(s|s)$ and $e(t)$ are independent. From \eqref{eq:5-201}, one can then compute $K(s,s|T)=\bE[\widehat{\epsilon}(s | T)\widehat{\epsilon}^{*}(s | T)]$ and obtain \eqref{eq:BF_formula} for the special case of $s=t$, as in \citet[Theorem 16.5.1]{Kailath2000}.

	\section{SUMMARY OF RESULTS}\label{sec:summary}
	In this section we summarize the findings and key expressions of the article. We also highlight the improvements made over previous formulas for kernels and controlled linear systems. Our problem was to find the operator $\hat{S}_s(t | T)$ minimizing the variance of the estimation error,
	\begin{equation}
		\epsilon_{S}(s | T)=x(s)-x_{S}(s|T)=\widetilde{x}(s)-\int_{t_0}^{T}S_s(t | T)d\widetilde{y}(t).\label{eq:2-9_summ}
	\end{equation}
	\begin{equation}
		\hat{S}_s(\cdot | T) \in \argmin_{S_s(\cdot | T)} \Gamma_{S}(s | T)=\bE[\epsilon_{S}(s | T)(\epsilon_{S}(s | T))^{*}].\label{eq:2-10_summ}
	\end{equation}
	The latter is related to the posterior covariance through the expression $\hat{S}_s(t|T)=K(s,t|T)H^{*}(t)R^{-1}(t)$ (Corollary~\ref{prop3-1}) with
	\begin{equation}
		K(s,t | T)=\bE[\epsilon_{\hat S}(s | T)(\epsilon_{\hat S}(t | T))^{*}]\in\cL(\R^{n,*},\R^n).\label{eq:kernel_def_summ}
	\end{equation}
	The explicit formula for $K$ can be obtained by finding the Green functions of a two-point boundary value problem:
	\begin{align}
		\dfrac{d\hat{\mu}}{dt}&=F(t)\hat{\mu}(t)-G(t)Q(t)G^{*}(t)\hat{\nu}(t)+l_\mu(t)\label{eq:twoPtsBdry_summ}\\
		-\dfrac{d\hat{\nu}}{dt}&=F^{*}(t)\hat{\nu}(t)+H^{*}(t)R^{-1}(t)H(t)\hat{\mu}(t)-l_\nu(t) \nonumber \\
		\hat{\mu}(t_0)&=-\Pi_{0}\hat{\nu}(t_0),\;\hat{\nu}(T)=\Sigma_T \hat{\mu}(T). \nonumber
	\end{align}
	The canonical resolution proceeds by introducing two Riccati differential equations
	\begin{align}
		-\dfrac{d}{dt}\Sigma&=\Sigma(t)F(t)+F^{*}(t)\Sigma(t)-\Sigma(t)G(t)Q(t)G^{*}(t)\Sigma(t)+H^{*}(t)R^{-1}(t)H(t), \quad &\Sigma(T)=\Sigma_T;\label{eq:2-20_summ}\\
		\dfrac{d}{dt}\Pi&=F(t)\Pi(t)+\Pi(t)F^{*}(t)-\Pi(t)H^{*}(t)R^{-1}(t)H(t)\Pi(t)+G(t)Q(t)G^{*}(t),\quad &\Pi(t_0)=\Pi_{0}.\label{eq:5-110_summ}
	\end{align}
	which we complement by looking for matrix-valued kernels $K$ and $\Lambda$ satisfying
	\begin{align}
		\hat{\mu}(s)=\int_{t_0}^{T}K(s,t | T)l_\nu(t)dt \quad &\text{ for } l_\mu(\cdot)\equiv 0, \label{eq:Wiener_K_summ}\\
		\hat{\nu}(s)=\int_{t_0}^{T}\Lambda(s,t | T)l_\mu(t)dt  \quad &\text{ for } l_\nu(\cdot)\equiv 0. \label{eq:Wiener_L_summ}
	\end{align}
	where $K$ corresponds to the covariance \eqref{eq:kernel_def_summ} as shown in Proposition~\ref{prop:covar_BF_formula}. Introducing a semigroup associated with the matrix function $F(t)-G(t)Q(t)G^{*}(t)\Sigma(t)$ denoted $\Phi_{F,\Sigma}(s,t)$ (resp.\ $F(s)-\Pi(s)H^{*}(s)R^{-1}(s)H(s)$, denoted $\Phi_{F,\Pi}(s,t)$), we obtained two symmetric formulas (Theorem~\ref{thm:main_noyaux}):
	\begin{align}
		K(s,t | T)&=\Phi_{F,\Sigma}(s,t_0)\Pi_{0}^{\frac{1}{2}}(\Id+\Pi_{0}^{\frac{1}{2}}\Sigma(t_0)\Pi_{0}^{\frac{1}{2}})^{-1}\Pi_{0}^{\frac{1}{2}}\Phi_{F,\Sigma}^{*}(t,t_0)+\int_{t_0}^{\min(s,t)}\Phi_{F,\Sigma}(s,\tau)G(\tau)Q(\tau)G^{*}(\tau)\Phi_{F,\Sigma}^{*}(t,\tau)d\tau\label{eq:3-3_summ}\\
		\Lambda(s,t|T)&=\Phi_{F,\Pi}^{*}(T,s)\Sigma_T^{\frac{1}{2}}(\Id+\Sigma_T^{\frac{1}{2}}\Pi(T)\Sigma_T^{\frac{1}{2}})^{-1}\Sigma_T^{\frac{1}{2}}\Phi_{F,\Pi}(T,t)+\int_{\max(s,t)}^{T}\Phi_{F,\Pi}^{*}(\tau,s)H^{*}(\tau)R^{-1}(\tau)H(\tau)\Phi_{F,\Pi}(\tau,t)d\tau\label{eq:5-16_summ}
	\end{align}
	For which we proved (Theorem~\ref{theo4-1}), that $K$ was the reproducing kernel of a space of controlled trajectories equipped with a quadratic norm
	\begin{multline}
		\Sx=\{x(\cdot)\in H^{1}(t_0,T;\R^{n})|\exists\:u(\cdot)\in L^{2}(t_0,T;\R^{p}),\xi\in \R^{n}\:\text{ s.t.\ }\\ \dfrac{d}{d\tau}x=F(\tau)x(\tau)+G(\tau)Q^{\frac{1}{2}}(\tau)u(\tau),x(t_0)=\Pi_{0}^{\frac{1}{2}}\xi\}\label{eq:4-1_summ}
	\end{multline}
	\begin{equation}
		\noridx{x(\cdot)}{\Sx}^2=J_x(\xi,u(\cdot))=\|\xi\|^{2}+\scalh{\Sigma_Tx(T)}{x(T)}+\int_{t_0}^{T}\|u(\tau)\|^{2}d\tau+\int_{t_0}^{T}\scalh{H^{*}(\tau)R^{-1}(\tau)H(\tau)x(\tau)}{x(\tau)}d\tau.\label{eq:3-19_summ}
	\end{equation}
	Similarly, $\Lambda$ is associated with a Hilbert space of information vectors (Theorem~\ref{theo5-2})
	\begin{multline}
		\Sla=\{\lambda(\cdot)\in H^{1}(t_0,T;\R^{n})|\:v(\cdot)\in L^{2}(t_0,T;\R^{m}),z\in\R^n\:\text{ s.t.\ }\\ -\dfrac{d}{dt}\lambda(t)=F^{*}(t)\lambda(t)+H^{*}(t)v(t),\;\lambda(T)=\Sigma_T^{\frac{1}{2}} z\}\label{eq:5-10_summ}
	\end{multline}
	\begin{equation}
		\noridx{\lambda(\cdot)}{\Sla}^2=J_\lambda(z,v(\cdot))=\scalh{\Pi_{0}\lambda(t_0)}{\lambda(t_0)}+\norh{z}^2+\int_{t_0}^{T}\scalh{G(t)Q(t)G^{*}(t)\lambda(t)}{\lambda(t)}dt+\int_{t_0}^{T}\scalh{R(t)v(t)}{v(t)}dt\label{eq:5-9_summ}
	\end{equation}
	The two kernels are related by the following formula (Proposition~\ref{prop:covar_BF_formula})
	\begin{equation}
		K(s,t | T)=\Pi(s)\Phi_{F,\Pi}^{*}(t,s)\1_{s\le t}+\Phi_{F,\Pi}(s,t)\Pi(t)\1_{s> t}-\Pi(s)\Lambda(s,t|T)\Pi(t).\label{eq:BF_formula_summ}
	\end{equation}
	For any realization $\widetilde{y}(\cdot)$ of $Y-\bar y$, the kernels allow to define two dual optimization problems (see Propositions~\ref{prop:primal_prob_rkhs} and \ref{prop:dual_prob_rkhs} for the definitions) formalizing \eqref{eq:mvle-ls}-\eqref{eq:mvle-3}
	\begin{align}
		&\min_{\tilde{x}(\cdot)\in\Sx} \noridx{R(t)^{\nicefrac{-1}{
						2}}\widetilde y(\cdot)}{L^2}^2+\noridx{\widetilde{x}(\cdot)}{\Sx}^2-2\scalidx{H^{*}(\cdot)R^{-1}(\cdot)\widetilde y(\cdot)}{\widetilde{x}(\cdot)}{L^2([t_0,T])} \label{eq:prob_primal_RKHS} \\ 
				&{\hspace{2cm}\scriptstyle  =\int_{t_{0}}^T \noridx{\tilde{y}(t)- H(t)\tilde{x}(t)}{R(t)^{-1}}^2dt+ \noridx{\tilde{x}(\cdot)}{\Sx}^2-\int_{t_{0}}^T \noridx{ H(t)\tilde{x}(t)}{R(t)^{-1}}^2dt}\nonumber \\
			 &\min_{\lambda(\cdot)\in\Sla}\noridx{\lambda(\cdot)}{\Sla}^2 - 2\int_{t_{0}}^T \scalh{R(t)\proj\nolimits^{\|\cdot\|_{R(t)}}_{\Img H(t)}(\widetilde{y}(t))}{v(t)}dt-\noridx{R(\cdot)^{\nicefrac{1}{
			 			2}}\proj\nolimits^{\|\cdot\|_{R(\cdot)}}_{\Ker H^*(\cdot)}(\widetilde{y}(\cdot))}{L^2}^2 \label{eq:prob_dual_RKHS}\\
			 &{\hspace{2cm}\scriptstyle =\int_{t_{0}}^T \noridx{\proj\nolimits^{\|\cdot\|_{R(t)}}_{\Img H(t)}(\tilde{y}(t))- v(t)}{R(t)}^2dt+\int_{t_0}^{T}\scalh{G(t)Q(t)G^{*}(t)\lambda(t)}{\lambda(t)}dt  +\scalh{\Pi_{0}\lambda(t_0)}{\lambda(t_0)}
			 	+\scalh{\Sigma_T^\ominus \lambda(T)}{\lambda(T)}-\noridx{R(\cdot)^{\nicefrac{-1}{
			 				2}}\tilde{y}(\cdot)}{L^2}^2}. \nonumber
	\end{align}
	This relates Kalman filtering to optimization problems over RKHSs, which act in filtering problems as the Sobolev spaces in calculus of variations. Inspired by Table 15.1 in \citet{Kailath2000}, which written for discrete-time estimation problems, we summarize in Table~\ref{table:four_prob} the relations between the deterministic optimal control problems written over RKHSs \eqref{eq:prob_primal_RKHS}-\eqref{eq:prob_dual_RKHS} and the original stochastic smoothing problems. The lower line, consisting of problems (iii) and (iv) over dual variables, is arguably less studied both in estimation and control, and we leave to future work the identification of the backward Markovian Gaussian process underlying $\Lambda$.
	
	\begin{table}[!h]
		\centering
		\resizebox{\columnwidth}{!}{%
		\begin{tabular}{|c|c|c|}
			\hline 
			& Stochastic problems & Deterministic problems \\
			\hline 
			\rotatebox[origin=r]{90}{ primal variables  } & \makecell[tl]{(i) Given Gaussian processes $(X_t)_{t\in[0,T]}$, $(Y_t)_{t\in[0,T]}$ \\Solve linear MMSE i.e.\ \eqref{eq:mvle}\\\\ $ \min_{\hat{X}\in \cL(Y)} \bE((X-\hat X)^\top(X-\hat X))$ \\\\ Optimum: $\hat{X}=\bE[X | Y]=\hat S Y$} & \makecell[tl]{(ii) Given RKHS $\Sx$ with kernel $K$, observations $\tilde{y}(\cdot)$ \\ Solve primal optimal control problem over trajectories \\\\ $ { \min_{\tilde{x}(\cdot)\in\Sx} \noridx{\tilde{x}(\cdot)}{\Sx}^2-2\scalidx{H^{*}(\cdot)R^{-1}(\cdot)\tilde y(\cdot)}{\tilde{x}(\cdot)}{L^2([t_0,T])}}$\\\\ Optimum: $\hat{x}(\cdot)=\int_{t_0}^{T}K(\cdot,t | T)H^{*}(t)R^{-1}(t)\widetilde{y}(t)dt$} \\
			\hline 
			\rotatebox[origin=r]{90}{ dual variables  } & \makecell[tl]{(iii) Given Gaussian processes $(X^\top_t)_{t\in[0,T]}$, $(Y_t)_{t\in[0,T]}$ \\ Solve over $(\Lambda_t)_{t\in[0,T]}$ \eqref{eq:mvle-dual2} \\\\ $ \min_{\Lambda \in \cL(Y)^\perp}\bE((X^\top-\Lambda)^\top(X^\top-\Lambda))$ \\\\ Optimum: $\hat{X}^\top=X^\top-\hat \Lambda$ } & \makecell[tl]{(iv) Given RKHS $\Sla$ with kernel $\Lambda$, observations $\tilde{y}(\cdot)$\\ Solve dual optimal control problem over adjoint/information \\\\ ${ \min_{\lambda(\cdot)\in\Sla}\noridx{\lambda(\cdot)}{\Sla}^2 - 2\int_{t_{0}}^T \scalh{R(t)\proj\nolimits^{\|\cdot\|_{R(t)}}_{\Img H(t)}(\tilde{y}(t))}{v(t)}dt}$\\ Optimum: $\hat v(t)=-R^{-1}(t)H(t)\hat x(t)+\proj\nolimits^{\|\cdot\|_{R(t)}}_{\Img H(t)}(\tilde{y}(t))$} \\ 
			\hline
		\end{tabular}
	}
		\caption{Summary of the four optimization problems considered in the article. One moves vertically by permuting min-max into max-min, the problems being (Fenchel) dual. One moves horizontally by formally setting $dw(t)=u(t)dt$ and considering the same kernel for (stochastic) covariance of optimal error and (deterministic) trajectories, the problems being (kernel) ``equivalent''.}\label{table:four_prob}
	\end{table}
	
	\vspace{3mm}
	\noindent\tb{Comparison with previous results on kernels and control:} The trajectory space seen as an RKHS was first presented in previous articles on linear-quadratic optimal control:
	\begin{itemize}
		\item \cite{aubin2020hard_control} first introduced the idea of considering the vector space $(\Sx,\scalidx{\cdot}{\cdot}{\Sx})$ as a reproducing
		kernel Hilbert space in the case $\Sigma_{T}=0$. A two-point boundary system to obtain $K$ was given but not solved explicitly. The closed form formula were only given for $\Pi_0=\Id$ and $H\equiv 0$. The kernel was used to guarantee the satisfaction of state constraints;
		\item \cite{aubin2020Riccati} observed that, for $\Sigma_{T}\neq 0$ and $\Pi_{0}\rightarrow\infty$, the map $t_0\mapsto K(t_0,t_0|[t_0,T])$ satisfied a forward Riccati equation and was the inverse of the usual backward Riccati matrix considered in linear-quadratic optimal control. It was underlined that kernels shift the focus on trajectories rather than their parametrization by controls;
		\item \cite{aubin2022operator} gave a closed form formula for $K$, when $\Sigma_{T}=0$, in the general case of an infinite-dimensional state with values in a Hilbert space to tackle linear PDE control. We emphasized there that this allowed for representer theorems and closed-form solutions when considering linear-quadratic optimal control problems.
	\end{itemize}
	To summarize, the kernels considered appear in linear-quadratic optimal control because of Hilbertian vector spaces of trajectories, while, for estimation problems, they appear through covariances of Gaussian processes. It is this ``dual'', deterministic and stochastic, nature of kernels which underlies the ``duality'' between optimal control and estimation in the Linear-Quadratic case.\\
	
	\noindent\tb{Conclusion.} We improved on our previous results by considering the dual RKHS $(\Sla,\scalidx{\cdot}{\cdot}{\Sla})$ of covectors/information vectors, by relaxing as much as possible the invertibility requirements on the matrices $Q,\Pi_0,\Sigma_{T}$, and most of all by considering an estimation, rather than optimal control, problem and consequently defining the kernel $K$ as the covariance of the estimation error. This way, we derived novel formulas for the covariances of the Markovian Gaussian processes induced by linear SDEs.  The kernels presented should also in principle allow for incorporating constraints or considering various sampling times of observations, and are computable through the matrix exponential formulas given in Appendix A.2. We considered here a continuous-time observation process, the extension to discrete-time measurements is straightforward by replacing integrals with finite sums when dealing with the observation operator $H$. Similarly generalization to an infinite-dimensional state space can be done in the spirit of \citet{aubin2022operator}. We could not discuss here the question of duality as a formal change of variables as in \citet{Kalman1961NewRI,Todorov2008}, as well as the relation between value function and likelihood. Indeed, the two Riccati equations outlined correspond respectively to the backward evolution of the Hessian of the quadratic value function through the Hamilton-Jacobi-Bellman equation, and to the forward evolution of the Hessian of the Gaussian posterior density through the Fokker-Planck equation. This paves the way to extending to nonlinear filtering the kernel viewpoint. We leave these interesting directions to future work.\\
	
	\noindent\tb{Acknowledgments:} We thank the anonymous referee for his positive and constructive comments. PCAF expresses his gratitude to Marc Lambert and Hans Kersting for the numerous discussions at the SIERRA laboratory on Kalman filtering, which spurred him into exploring the duality and Bayesian aspects.
	\section*{Appendix}
	\subsection*{A.1 Existence of the solution of the Riccati equation}
	\begin{lem} \label{lem3-1} The Riccati equation (\ref{eq:2-20}) has a solution
		$\Sigma(t)$ which is symmetric and positive semi-definite on the interval $[t_0,T]$. 
	\end{lem}
	\begin{proof}
		The differential equation \eqref{eq:2-20} has a unique local solution on an interval
		$[t_{1},T]$ for $t_{1}$ sufficiently close to $T$. The solution
		is symmetric, since the transpose satisfies the equation. We can then
		consider the differential equation 
		\begin{equation}
			\dfrac{d\gamma}{dt}=(F(t)-G(t)Q(t)G^{*}(t)\Sigma(t))\gamma(t),\quad \gamma(t_{1})=\gamma_{0}.\label{eq:2-21}
		\end{equation}
		We can then compute 
		\begin{align*}
			\dfrac{d\scalh{\gamma(t)}{\Sigma(t)\gamma(t)}}{dt}&=2\scalh{\gamma(t)}{\Sigma(t)\dfrac{d\gamma(t)}{dt}}+\scalh{\gamma(t)}{\dfrac{d\Sigma(t)}{dt}\gamma(t)}\\
			&=-\scalh{\gamma(t)}{\left(\Sigma(t)G(t)Q(t)G^{*}(t)\Sigma(t)+H^{*}(t)R^{-1}(t)H(t)\right)\gamma(t)}\leq0
		\end{align*}
		and since $\Sigma(T)\succcurlyeq 0$, by integration between $t_{1}$ and $T$, we
		obtain $\scalh{\gamma_0}{\Sigma(t_1)\gamma_0}\geq0$. Since $\gamma_{0}$
		is arbitrary, we obtain $\Sigma(t_{1})\geq0$. We could have started in
		(\ref{eq:2-21}) at any point on the interval $(t_{1},T)$. Therefore
		the local solution satisfies $\Sigma(t)\geq0$. Consider next the semi group $\Phi_F(t,t_{1})$,
		on the interval $t\in(t_{1},T)$ defined by the differential equation 
		\begin{equation}
			\dfrac{d}{dt}\Phi_F(t,t_{1})=F(t)\Phi_F(t,t_{1}), \quad \Phi_F(t_{1},t_{1})=\Id.\label{eq:2-22}
		\end{equation}
		If we consider 
		\[
		\dfrac{d\scalh{\Phi_F(t,t_{1})\gamma_{0}}{\Sigma(t)\Phi_F(t,t_{1})\gamma_{0}}}{dt}=\scalh{\Phi_F(t,t_{1})\gamma_{0}}{\left(F^{*}(t)\Sigma(t)+\Sigma(t)F(t)+\dfrac{d}{dt}\Sigma(t)\right)\Phi_F(t,t_{1})\gamma_{0}},
		\]
		from the Riccati equation, we obtain
		\[
		\dfrac{d\scalh{\Phi_F(t,t_{1})\gamma_{0}}{\Sigma(t)\Phi_F(t,t_{1})\gamma_{0}}}{dt}\geq-\scalh{\Phi_F(t,t_{1})\gamma_{0}}{H^{*}(t)R^{-1}(t)H(t)\Phi_F(t,t_{1})\gamma_{0}}
		\]
		Integrating between $t_{1}$ and $T$, we get 
		\begin{equation}
			\scalh{\gamma_{0}}{\Sigma(t)\gamma_{0}}\leq\int_{t_{1}}^{T}\scalh{\Phi_F(t,t_{1})\gamma_{0}}{H^{*}(t)R^{-1}(t)H(t)\Phi_F(t,t_{1})\gamma_{0}}\leq C\|\gamma_{0}\|^{2}\label{eq:2-23}
		\end{equation}
		for some constant $C>0$. It follows that the solution $\Sigma(t)$ can be extended beyond $t_{1}$
		and finally up to $t_0$. This completes the proof. 
	\end{proof}

	\subsection*{A.2 Computation of kernels of linear SDEs through exponentials of the Hamiltonian matrix}
	
	The following method to compute Gramians can be traced at least back to \citet{VanLoan1978computingExpM} in the time-invariant case and was mentioned in \citet[Section 8]{kalman1960contribution} as a periodically rediscovered way to solve differential Riccati equations. It made its way in the linear SDE literature \citep[see e.g.][Section 6.3, p.84]{Srkk2019} where $H\equiv 0$ and $\Sigma \equiv 0$. However a more systematic presentation through the Hamiltonian matrix can be found in \citet{speyer2010primer} for linear-quadratic control. In particular this allows to consider the important case of $GQG^*\not\equiv 0$ and $H\not\equiv 0$. We focus below on computing $K$ but similar operations can be performed to obtain $\Lambda$ by changing variables, as was done in the proof of Theorem~\ref{theo5-2}. Consider the equations
	\begin{align*}
		\dfrac{d\hat{\mu}}{dt}&=F(t)\hat{\mu}(t)-G(t)Q(t)G^{*}(t)\hat{\nu}(t), \quad &&\hat{\mu}(t_0)=- \Pi_0 \hat{\nu}(t_0) \\
		\dfrac{d\hat{\nu}}{dt}&=-F^{*}(t)\hat{\nu}(t)-H^{*}(t)R^{-1}(t)H(t)\hat{\mu}(t), \quad &&\hat{\nu}(T)=\Sigma_T \hat{\mu}(T)
	\end{align*}
	and write them in matrix form introducing the Hamiltonian matrix $\bH(t)$
	\begin{equation}
		\dfrac{d}{dt}\begin{pmatrix}
			\hat{\mu}(t)\\
			\hat{\nu}(t)
		\end{pmatrix}=\underbrace{\begin{pmatrix} F(t) &  -G(t)Q(t)G^{*}(t)\\ -H^{*}(t)R^{-1}(t)H(t) & -F^{*}(t) \end{pmatrix}}_{\bH(t)}\begin{pmatrix}
			\hat{\mu}(t) \\
			\hat{\nu}(t)
		\end{pmatrix}.
	\end{equation} 
	Denote by $\Phi_{\bH}(T,t)$ the transition matrix  $\frac{d }{dt}\Phi_{\bH}(T,t)=-\bH(t)\Phi_{\bH}(T,t)$, $\Phi_{\bH}(T,T)=\Id$ and set as in \citet[eq.(5.72)]{speyer2010primer}
	$$\Phi_{\bH}(T,t):= \begin{pmatrix} \Phi_{11}(T,t) &  \Phi_{12}(T,t)\\ 
		\Phi_{21}(T,t) & \Phi_{22}(T,t) \end{pmatrix} \quad \bar \Phi_{\bH}(T,t):= \begin{pmatrix} \bar\Phi_{11}(T,t) &  \bar\Phi_{12}(T,t)\\ 
		\bar	\Phi_{21}(T,t) &\bar \Phi_{22}(T,t) \end{pmatrix}= \begin{pmatrix} \Id & 0\\ 
		-\Sigma_T & \Id \end{pmatrix} \cdot \Phi_{\bH}(T,t).$$
	Assume that $\bar \Phi_{22}(T,t)$ is invertible, then through some calculations one can show that \citep[see][Chapter 5, eq.(5.77,5.160, 5.163)]{speyer2010primer}
	\begin{gather*}
		\Sigma(t)=-\bar\Phi_{22}^{-1}(T,t)\bar\Phi_{21}(T,t),\\
		\Phi_{F,\Sigma}(t,T)=\bar \Phi_{22}^{*}(T,t),\\
		\int_{T}^{s}\Phi_{F,\Sigma}(T,\tau)G(\tau)Q(\tau)G^{*}(\tau)\Phi_{F,\Sigma}^{*}(T,\tau)d\tau=\bar\Phi_{12}(T,s)\bar\Phi_{22}^{-1}(T,s).
	\end{gather*}
	where $\Sigma(t)$ is the solution of the backward Riccati equation \eqref{eq:2-20_summ} and $\Phi_{F,\Sigma}(s,t)$ the semigroup associated with $F(t)-G(t)Q(t)G^{*}(t)\Sigma(t)$. Define two auxiliary kernels which we will compute independently.
	\begin{align}
		K(s,t | T)&=\underbrace{\Phi_{F,\Sigma}(s,t_0)\Pi_{0}^{\frac{1}{2}}(\Id+\Pi_{0}^{\frac{1}{2}}\Sigma(t_0)\Pi_{0}^{\frac{1}{2}})^{-1}\Pi_{0}^{\frac{1}{2}}\Phi_{F,\Sigma}^{*}(t,t_0)}_{=:K_0(s,t)}+\underbrace{\int_{t_0}^{\min(s,t)}\Phi_{F,\Sigma}(s,\tau)G(\tau)Q(\tau)G^{*}(\tau)\Phi_{F,\Sigma}^{*}(t,\tau)d\tau}_{=:K_1(s,t)}\label{eq:def_kernelAux}
	\end{align}
	By the above expressions, we have for $s\le t$
	\begin{align*}
		K_{0}(s,t)&=\bar\Phi_{22}^{*}(T,s)\bar\Phi_{22}^{-1,*}(T,t_0)\Pi_{0}^{\frac{1}{2}}(\Id+\Pi_{0}^{\frac{1}{2}}\bar\Phi_{22}^{-1}(T,t_0)\bar\Phi_{21}(T,t_0)\Pi_{0}^{\frac{1}{2}})^{-1}\Pi_{0}^{\frac{1}{2}}\bar\Phi_{22}^{-1}(T,t_0)\bar\Phi_{22}(T,t)\\
		K_{1}(s,t)&=\bar\Phi_{22}^{*}(T,s)[\bar\Phi_{12}(T,s)\bar\Phi_{22}^{-1}(T,s)-\bar\Phi_{12}(T,t_{0})\bar\Phi_{22}^{-1}(T,t_{0})]\bar\Phi_{22}(T,t)
	\end{align*}
	Formulas for $s>t$ can be obtained by Hermitian symmetry, $K(s,t|T)=K(t,s|T)^*$. These formulas simplify drastically when $H(\cdot)\equiv0$ and $\Sigma_T=0$, so $\Sigma(\cdot)\equiv 0$, $\bar \Phi_{22}^{*}(T,t)=\Phi_{F}(t,T)$
	(alternative computations could be done in this case with $\Phi_{\bH}(t,t_0)$) which inverse is easily computable. Otherwise every inversion should be computed numerically. In the time-invariant case, we have $\Phi_{\bH}(T,t)=e^{(T-t)\bH}$ and all quantities can be computed through matrix exponentials.
	
	\bibliographystyle{unsrtnat}
	\bibliography{biblioPCAF,biblioKalmanFilter} 
	
\end{document}